\documentclass{article}
\usepackage{setspace}
\usepackage[utf8]{inputenc}
\usepackage[english]{babel}
\usepackage{amsfonts,amsmath,amssymb,mathtools}
\usepackage[left=2cm,right=2cm,top=2cm,bottom=2cm,bindingoffset=0cm]{geometry}
\usepackage{graphicx}
\usepackage{authblk}
\usepackage{xcolor}
\usepackage{hyperref}
\usepackage{dsfont}
\usepackage{authblk}
\usepackage{slashed}
\usepackage{braket}
\usepackage{comment}
\usepackage{pb-diagram}
\usepackage{wasysym}
\usepackage{amsthm}
\usepackage{cancel}
\usepackage{appendix}
\usepackage{tikz}\usetikzlibrary{tikzmark}\usetikzlibrary{calc}
\usetikzlibrary{fadings} 
\usetikzlibrary{positioning}% To get more advances positioning options
\usetikzlibrary{backgrounds} % LATEX and plain TEX
\usetikzlibrary{decorations.pathreplacing}
\usetikzlibrary{arrows}% To get more arrow heads
\definecolor{airforceblue}{rgb}{0.36, 0.54, 0.66}	
\definecolor{beige}{rgb}{0.96, 0.96, 0.86}
\definecolor{bittersweet}{rgb}{1.0, 0.44, 0.37}
\definecolor{melon}{rgb}{0.99, 0.74, 0.71}
\definecolor{mustard}{rgb}{1.0, 0.86, 0.35}
\definecolor{lava}{rgb}{0.81, 0.06, 0.13}
\definecolor{magnolia}{rgb}{0.97, 0.96, 1.0}
\definecolor{lavendermist}{rgb}{0.9, 0.9, 0.98}
\definecolor{lavendergray}{rgb}{0.77, 0.76, 0.82}
\definecolor{palepink}{rgb}{0.98, 0.85, 0.87}
\definecolor{palesilver}{rgb}{0.79, 0.75, 0.73}
\definecolor{cadetgrey}{rgb}{0.57, 0.64, 0.69}
\definecolor{anti-flashwhite}{rgb}{0.95, 0.95, 0.96}
\colorlet{Light0anti-flashwhite}{anti-flashwhite!70!white}
\colorlet{Lightanti-flashwhite}{anti-flashwhite!50!white}
\colorlet{Light2anti-flashwhite}{anti-flashwhite!30!white}
\definecolor{linkcolor}{rgb}{0,0,1}
\definecolor{urlcolor}{rgb}{0,0,1}

\usepackage{tikz-cd}
\usetikzlibrary{decorations.markings}
\usetikzlibrary{positioning}
\usepackage{braids}
\usepackage{array}
\usetikzlibrary{calc, knots, hobby, arrows.meta, through}
\pgfdeclarelayer{background layer}
\pgfdeclarelayer{foreground layer}
\pgfsetlayers{background layer,main,foreground layer}
\tikzstyle{arrow} = [thick,-{Stealth[length=5mm]}]
\usepackage{ytableau}
\usepackage{longtable}
\usepackage{empheq}
\usepackage{multicol}
\usepackage{mathrsfs}
\usepackage{diagbox}
\usepackage[vcentermath]{youngtab}
\usepackage[natbib=true, backend = bibtex, style=numeric-comp, sorting=none]{biblatex}
\hypersetup{pdfstartview=FitH,  linkcolor=linkcolor,urlcolor=urlcolor,citecolor=urlcolor, colorlinks=true}

\newcommand\bem{\begin{pmatrix}}
\newcommand\eem{\end{pmatrix}}
\newcommand\beq{\begin{equation}}
\newcommand\eeq{\end{equation}}
\newcommand\beqs{\begin{equation*}}
\newcommand\eeqs{\end{equation*}}

\date{}

\def\be{\begin{eqnarray}}
\def\ee{\end{eqnarray}}

\def\Tr{{\rm Tr}\,}

\definecolor{red}{rgb}{1,0,0}
\definecolor{orange}{rgb}{1,0.5,0}
\definecolor{violet}{rgb}{0.7,0,1}

\newtheorem{definition}{Definition}[section]

\newtheorem{remark}{Remark}[section]

\newcommand{\bbl}{{{-\hspace{-1.5mm}\bigcirc\hspace{-1.5mm}-}}}

\newcommand{\bubble}[2]{
	\coordinate (a123) at #1;
	\coordinate (b123) at #2;
	\coordinate (aa123) at ($(a123)!0.3!(b123)$);
	\coordinate (bb123) at ($(a123)!0.7!(b123)$);
	\node [draw] at ($0.5*(a123) + 0.5*(b123)$) [circle through={(aa123)}] {};
	\draw (a123) edge (aa123);
	\draw (b123) edge (bb123);
}

\newcommand{\xVogel}[4]{
	\coordinate (aaaa1) at #1;
	\coordinate (bbbb1) at #2;
	\coordinate (cccc1) at #3;
	\coordinate (mmmm1) at ($0.7*(cccc1) + 0.15*(aaaa1)+0.15*(bbbb1)$);
	\draw (aaaa1) -- (mmmm1) -- (cccc1);
	\draw (bbbb1) -- (mmmm1);
	\coordinate (eeee1) at ($0.4*(aaaa1)+0.6*(mmmm1)$);
	\coordinate (ffff1) at ($0.4*(bbbb1)+0.6*(mmmm1)$);
	\draw (eeee1)--(ffff1);
	\coordinate (gggg1) at ($0.65*(aaaa1)+0.35*(mmmm1)$);
	\coordinate (hhhh1) at ($0.65*(bbbb1)+0.35*(mmmm1)$);
	\draw (gggg1) -- (hhhh1);
	\foreach \x [evaluate=\x]  in {1-1, ...-1, #4-1} 
	\draw (${\x*1/(#4-1)}*(eeee1)+{1-\x*1/(#4-1)}*(ffff1)$) -- (${\x*1/(#4-1)}*(gggg1)+{1-\x*1/(#4-1)}*(hhhh1)$); 
}

\begin{document}

\title{\bf Construction of Lie algebra weight system kernel \newline via Vogel algebra}

\author[1,2]{{\bf D. Khudoteplov}\thanks{\href{mailto:khudoteplov.dn@phystech.edu}{khudoteplov.dn@phystech.edu}}}
\author[1,2,3]{{\bf E. Lanina}\thanks{\href{mailto:lanina.en@phystech.edu}{lanina.en@phystech.edu}}}
\author[1,2,3]{{\bf A. Sleptsov}\thanks{\href{mailto:sleptsov@itep.ru}{sleptsov@itep.ru}}}

\vspace{4.5cm}

\affil[1]{Moscow Institute of Physics and Technology, 141700, Dolgoprudny, Russia}
\affil[2]{Institute for Information Transmission Problems, 127051, Moscow, Russia}
\affil[3]{NRC "Kurchatov Institute", 123182, Moscow, Russia\footnote{former Institute for Theoretical and Experimental Physics, 117218, Moscow, Russia}}
\renewcommand\Affilfont{\itshape\small}

\maketitle

\vspace{-7cm}

\begin{center}
	\hfill MIPT/TH-27/24\\
	\hfill ITEP/TH-34/24\\
	\hfill IITP/TH-29/24
\end{center}

\vspace{4.0cm}

\begin{abstract}

{We develop a method of constructing a kernel of Lie algebra weight system. A main tool we use in the analysis is Vogel's $\Lambda$ algebra and the surrounding framework. As an example of a developed technique we explicitly provide all Jacobi diagrams lying in the kernel of $\mathfrak{sl}_N$ weight system at low orders. We also discuss consequences of the presence of the kernel in Lie algebra weight systems for detection of correlators in the 3D Chern--Simons topological field theory and for distinguishing of knots by the corresponding quantum knot invariants. 
}
\end{abstract}

\tableofcontents

\section{Introduction}

In 1984, V. Jones discovered a new polynomial knot invariant~\cite{jones2005jones,jones1987hecke,jones1985polynomial}, the Jones polynomial, associated with the Lie algebra $\mathfrak{sl}_2$ and its standard (fundamental) representation. Soon, generalizations of the Jones polynomials appeared: first to the algebras $\mathfrak{sl}_n$ and $\mathfrak{so}_n$ (the HOMFLY~\cite{freyd1985new,przytycki1987kobe} and Kauffman polynomials~\cite{kauffman1990invariant}, respectively), and then to arbitrary finite-dimensional representations by methods of quantum groups and $\cal R$-matrices (Reshetikhin--Turaev invariants~\cite{TURAEV1992865,Reshetikhin:1990pr}). The Reshetikhin--Turaev construction has the property of universality, since it is formulated in fairly general terms and defines knot invariants related to simple Lie algebras and their finite-dimensional representations. An extension of the construction to Lie superalgebras is also allowed \cite{zhang1991universal, rozansky1993s,geer2005kontsevich,viro2007quantum,queffelec2019note}.

In parallel, physicists A. Schwartz~\cite{schwarz1978partition, schwarz1979partition, schwarz1987baku} and E. Witten~\cite{Witten1988hf} discovered that the 3-dimensional quantum Chern--Simons field theory also allows one to construct these invariants. In this theory it is necessary to construct a gauge-invariant quantity called the vacuum expectation value of the Wilson loop, well known to physicists from the theory of strong interactions (i.e. quantum chromodynamics). This approach also has the property of universality and, on the one hand, opens the way for physicists to the world of low-dimensional topology, and, on the other hand, enriches the mathematical tools with new powerful methods coming from quantum field theory and string theory~\cite{witten1988topological,Witten1988hf,Chern1974ft,Guadagnini:1989kr,GUADAGNINI1990575,Kaul1991np,RamaDevi1992np,Ramadevi1993np,Ramadevi1994zb,Ooguri1999bv, Labastida2000zp, Labastida2000yw, Labastida2001ts, Marino2001re}. These invariants are now called in literature either Witten--Reshetikhin--Turaev invariants or quantum knot invariants. All these methods and approaches to constructing knot invariants, inspired by the ideas of quantum mechanics and quantum field theory and applied to problems of low-dimensional topology, formed a new direction called quantum topology.

In 1992, Birman and Lin discovered that the coefficients of a certain perturbative expansion of the HOMFLY polynomials are Vassiliev invariants~\cite{birman1993knot}, which V.Vassiliev defined in 1990 studying complements of discriminants in spaces of maps ~\cite{vassiliev1990cohomology} (more about Vassiliev invariants in the next section). In 1995, D. Bar-Nathan proved that for any quantum knot invariant the expansion coefficient is a Vassiliev invariant~\cite{bar1995vassiliev}. The inverse question naturally arises: Is any Vassiliev invariant some coefficient (or a linear combination of them) of the perturbative expansion of a suitable quantum invariant (or even a linear combination of them)?

This question can be formulated more precisely if we introduce the concept of the Kontsevich integral, which is a universal Vassiliev invariant~\cite{kontsevich1993vassiliev}, see Section~\ref{sec:VK-theorem}. Roughly speaking, this is a generating series for all Vassiliev invariants~\eqref{ki}. The series goes over the formal variable $\hbar$. Each Vassiliev invariant is weighted by a certain chord diagram, since the algebra of Vassiliev invariants is isomorphic to the algebra of chord diagrams (more in Section \ref{sec:prereq}). Any finite-dimensional Lie algebra (with a non-degenerate Killing form) defines a linear function (called weight system) on the algebra of chord diagrams. If we take $\mathfrak{sl}_n$ as such an algebra, then the Kontsevich integral turns into the HOMFLY polynomial, into the Kauffman polynomial in the case of $\mathfrak{so}_n$, etc. More precisely, it turns into a series of perturbative decomposition for this polynomial, see Section~\ref{sec:CS-KI}. The question is which chord diagrams lie in the kernel of this weight system. Vassiliev invariants that correspond to these diagrams will be impossible to extract from the perturbative decomposition of the original quantum invariant. Next, we can ask whether there are such chord diagrams that lie in the kernel of all weight systems constructed from finite-dimensional Lie algebras and superalgebras.

It turns out that such diagrams exist, what was proved in a seminal paper by Pierre Vogel \cite{vogel2011algebraic} in 2011. However, Vogel did much more than just prove this fact. He developed a diagrammatic technique for explicitly constructing such kernels for any finite dimensional Lie (super)algebra. In this paper we apply his technique to construct the kernel of the $\mathfrak{sl}_n$ weight system. Without Vogel's technique this problem is quite difficult, because it is difficult to calculate a Lie algebra weight system of all basis diagrams of a given order. Basically there are two ways: the \textit{universal} weight system and the weight system associated with a \textit{representation}. 

For the universal way one parameterizes the center of universal enveloping algebra for the corresponding Lie algebra and expresses the values of the weight system for all (basis) chord diagrams in terms of selected parameters. Then, it is easy to find the dimension of the weight system and therefore its kernel. In this way all calculations must be done analytically. It was suggested by D. Bar-Natan in \cite{bar1995vassiliev}, a specific calculation method was also given there. Further it was developed by Zh. Yang in~\cite{yang2023new} and M. Kazarian and S. Lando in~\cite{kazaryan2023weight}. However, as far as we know, the 8-th degree is already unattainable for this method at the current progress, but it is at the 8-th order that the first primitive element of the $\mathfrak{sl}_n$ weight system kernel appears.

Another way to compute a Lie algebra weight system is to provide it with a finite-dimensional representation. Then, the computation need not be done analytically, which is an advantage over the previous method. However, it also introduces a disadvantage, since it is necessary to ensure that the relations found are independent on the representation. As far as we know, this method has not been developed much. Some computations can be found in Alvarez's and Labastida's paper~\cite{alvarez1995numerical} for the fundamental representation of $\mathfrak{sl}_n$ and $\mathfrak{so}_n$, and for all representations in $\mathfrak{sl}_2$.

An interesting development is a paper \cite{lanina2022implications}, where a combination of both methods is used for the calculation, and in addition various symmetry considerations are used for the colored HOMFLY polynomials (i.e. $\mathfrak{sl}_n$ case), which impose certain restrictions on the possible form of weight systems. The developed method allows to find a large number of relations for the weight system $\mathfrak{sl}_n$ and to obtain a good upper bound on the kernel dimension, which is presented below in Table \ref{tab:dimensions}.

It is also important to note that the first person to use Vogel's technique to construct the kernel of weight systems was Lieberum~\cite{lieberum1999vassiliev}. However, he studied those elements that lie in the kernel of all Lie algebras at once, except $D(2,1,\alpha)$.

The simplest primitive Jacobi diagram (see Section \ref{sec:prereq} for the definitions), which belongs to the kernel of all weight systems of Lie algebras, found by P. Vogel, is of the 17-th order. Lieberum's simplest diagram is of the 15-th order. Today, these are too high orders to work with. For example, it is very difficult to understand whether there are relations for any two or three diagrams constructed by Vogel's method. To answer this question it is necessary to use AS and IHX relations, which is difficult. It is also not yet clear whether it is possible to find two such knots that do not differ in any invariants except those that correspond to the found diagrams. Moreover, it is not entirely clear how to approach this problem. When we restrict ourselves to the $\mathfrak{sl}_n$ weight system, all these problems are simplified, but not trivialized, which allows us to develop and refine the methods. First, the kernel of primitive diagrams  appears at the 8-th order, which is much simpler than the 17-th. Second, the weight system of $\mathfrak{so}_n$ or other algebras can be used to detect elements of the $\mathfrak{sl}_n$ kernel. Third, there are several candidates among the 13-crossing knots for the role of those that are not detected by some colored quantum $\mathfrak{sl}_n$ invariants and are not mutants. At least these knots have the same colored HOMFLY polynomials (those that we can calculate) and the uncolored Kauffman polynomial. Unfortunately, we do not yet know how to calculate the colored Kauffman polynomials well enough to test this hypothesis, but at least this task is quite straightforward. In addition, the 8-th order of the diagram allows us to expect that such knots, if they exist at all, will be encountered much earlier (in the sense that their minimum number of crossings should be smaller) than knots for the diagram of order 15 and especially of order 17.

This paper is organized as follows. In Section~\ref{sec:prereq}, main provisions of the theory of finite order invariants from the definition of Vassiliev invariants to Jacobi diagrams and weight systems are briefly outlined. Section~\ref{sec:Vogel-algebra} highlights results obtained by Pierre Vogel~\cite{vogel2011algebraic,vogel1999universal}. The definition of the algebra $\Lambda$ and the characters of Lie algebras on $\Lambda$ are given. Particular attention is paid to Vogel parameterization and the related concept of Vogel universality. In Section~\ref{sec:weight-system-kernel}, Vogel's results are used to develop a method for construction of Lie algebra weight
systems kernel. A method for constructing Jacobi diagrams that are vanished by weight systems is described. As an example we provide explicit Jacobi diagrams forming the kernel of $\mathfrak{sl}_n$ weight system at low orders. For diagrams constructed in this way, their linear independence is proven. Section~\ref{sec:CS-KI} is devoted to connection of Vassiliev invariants, chord diagrams, the Kontsevich integral and the results of this paper with the 3-dimensional topological Chern--Simons theory. In Section~\ref{sec:knot-det}, we discuss that the existence of Lie algebras weight systems kernels implies that there exist linear combination of knots that are distinguished by Vassiliev invariants coming from kernel diagrams but not by other Vassiliev invariants of the same or lower orders. However, there is still no argument that higher order Vassiliev invariants, and thus the corresponding quantum knot invariants, do not distinguish this linear combination of knots.

\setcounter{equation}{0}
\section{Vassiliev invariants and related notions}\label{sec:prereq}

In this section we introduce main objects of our study such as Vassiliev invariants, chord diagrams, Jacobi diagrams and weight systems. We mostly follow \cite{chmutov2012introduction}. 
\subsection{Vassiliev invariants}
	\begin{definition}
		A knot is an isotopy class of $S^1$ embeddings in $\mathbb{R}^3$.
	\end{definition}
	
	\begin{definition}
		A knot invariant is a function on the set of knots.
	\end{definition}
	
	\noindent Before defining Vassiliev invariants, it is necessary to expand the concept of a knot to include singular knots. 
	
	\begin{definition}
		A singular knot is an isotopy class of $S^1$ immersions in $\mathbb{R}^3$ such that all self-intersection points are simple double points with transversal intersections.
	\end{definition}

	\noindent The continuation of knot invariants to the set of singular knots is determined by \textit{Vassiliev skein relation} shown in Fig.\,\ref{fig:Vass-skein}. The Vassiliev skein relation is a rule for resolving double points of a singular knot. It can be easily seen that a value of an invariant at a singular knot does not depend on an order in which its double points are resolved.
	\begin{figure}[h]
		\centering
		\begin{tikzpicture}[scale=0.4]
			\coordinate (a) at (-7, 0);
			\coordinate (b) at (-0.5, 0);
			\coordinate (c) at (6, 0);
			\begin{pgfonlayer}{background layer}
				\draw [-{Stealth[length=3mm]}, ultra thick] ($(c)-(1, 1)$) -- ($(c)+(1, 1)$);
				\draw [-{Stealth[length=3mm]}, ultra thick] ($(b)-(-1, 1)$) -- ($(b)+(-1, 1)$);
				\draw [-{Stealth[length=3mm]}, ultra thick] ($(a)-(1, 1)$) -- ($(a)+(1, 1)$);
				\draw [-{Stealth[length=3mm]}, ultra thick] ($(a)-(-1, 1)$) -- ($(a)+(-1, 1)$);
			\end{pgfonlayer}
			\begin{pgfonlayer}{main}
				\draw [white, line width=10pt] ($(c)+(1, -1)$) -- ($(c)+(-1, 1)$);
				\draw [-{Stealth[length=3mm]}, ultra thick] ($(c)+(1, -1)$) -- ($(c)+(-1, 1)$);
				
				\draw [white, line width=10pt] ($(b)+(-1, -1)$) -- ($(b)+(1, 1)$);
				\draw [-{Stealth[length=3mm]}, ultra thick] ($(b)+(-1, -1)$) -- ($(b)+(1, 1)$);
				\draw [fill] (a) circle(0.15);	
			\end{pgfonlayer}
			\begin{pgfonlayer}{foreground layer}
				\draw [dashed] (a) circle (1.4142);
				\draw [dashed] (b) circle (1.4142);
				\draw [dashed] (c) circle (1.4142);
				\node [scale =1.6] at ($(a)+(-2.4, 0)$) { \textit{v} } ;
				\node [scale =2.5] at ($(a)+(-1.7, 0)$) { ( } ;
				\node [scale =2.5] at ($(a)+(1.7, 0)$) { ) } ;
				\node [scale =1.6] at ($(b)+(-2.4, 0)$) { \textit{v} } ;
				\node [scale =2.5] at ($(b)+(-1.7, 0)$) { ( } ;
				\node [scale =2.14] at ($(b)+(1.7, 0)$) { ) } ;
				\node [scale =1.6] at ($(c)+(-2.4, 0)$) { \textit{v} } ;
				\node [scale =2.5] at ($(c)+(-1.7, 0)$) { ( } ;
				\node [scale =2.5] at ($(c)+(1.7, 0)$) { ) } ;
				\node [scale=1.8] at ($0.6*(b) + 0.4*(c)+(0.2, 0)$) {$-$};
				\node [scale=1.8] at ($0.6*(a) + 0.4*(b)+(0.2, -0.05)$) {$=$};
			\end{pgfonlayer}
		\end{tikzpicture}
\vspace{-3mm}
		\caption{Vassiliev skein relation.}
		\label{fig:Vass-skein}
	\end{figure}

	\begin{definition}[\textbf{Vassiliev invariants}\label{def:Vass_inv_def} \cite{vassiliev1990cohomology}]
		An invariant $v$ is called a Vassiliev invariant of order no more than $n$ if $v$ vanishes at singular knots with more than $n$ double points. The space of Vassiliev invariants of order $\leq n$ is denoted $\mathcal{V}_n$.
	\end{definition}

	\noindent Vassiliev invariants form a filtered vector space, and they also form a bialgebra (see \cite[\S 4.3]{chmutov2012introduction}).

	\subsection{Chord diagrams}
\label{sec:chord_diag}
	
	\begin{definition}
		A chord diagram of order $n$ is an oriented circle with $n$ pairs of distinct points. The set of chord diagrams of order $n$ is denoted $\mathbf{A}_n$.
	\end{definition}

	Chord diagrams naturally arise in the theory of Vassiliev invariants \cite{vassiliev1990cohomology}. To demonstrate this relationship, one needs to consider $v \in \mathcal{V}_n$. Suppose $K$ is an arbitrary singular knot with exactly $n$ double points. Next, let $K'$ be a knot that is identical to $K$ in all regions of space except for one where they differ from each other in the same way as the two diagrams in Fig.\,\ref{fig:KK'} do. We also define $K^\bullet$ to be a knot with a self-intersection in the area where $K$ and $K'$ differ, and identical to both in the remaining space (see Fig.\,\ref{fig:KK'}).

	\begin{figure}[!h]
		\centering
		\begin{tikzpicture}[scale=0.4]
			\coordinate (a) at (5, 0); %K^{\bullet}
			\coordinate (b) at (-5, 0); % K
			\coordinate (c) at (0, 0); % K'
			\begin{pgfonlayer}{background layer}
				\draw [-{Stealth[length=3mm]}, ultra thick] ($(c)-(1, 1)$) -- ($(c)+(1, 1)$);
				\draw [-{Stealth[length=3mm]}, ultra thick] ($(b)-(-1, 1)$) -- ($(b)+(-1, 1)$);
				\draw [-{Stealth[length=3mm]}, ultra thick] ($(a)-(1, 1)$) -- ($(a)+(1, 1)$);
				\draw [-{Stealth[length=3mm]}, ultra thick] ($(a)-(-1, 1)$) -- ($(a)+(-1, 1)$);
			\end{pgfonlayer}
			\begin{pgfonlayer}{main}
				\draw [white, line width=10pt] ($(c)+(1, -1)$) -- ($(c)+(-1, 1)$);
				\draw [-{Stealth[length=3mm]}, ultra thick] ($(c)+(1, -1)$) -- ($(c)+(-1, 1)$);
				
				\draw [white, line width=10pt] ($(b)+(-1, -1)$) -- ($(b)+(1, 1)$);
				\draw [-{Stealth[length=3mm]}, ultra thick] ($(b)+(-1, -1)$) -- ($(b)+(1, 1)$);
				\draw [fill] (a) circle(0.15);	
			\end{pgfonlayer}
			\begin{pgfonlayer}{foreground layer}
				\draw [dashed] (a) circle (1.4142);
				\draw [dashed] (b) circle (1.4142);
				\draw [dashed] (c) circle (1.4142);
				\node at ($(b)+(0, -2)$) {$K$};
				\node at ($(c)+(0, -2)$) {$K'$};
				\node at ($(a)+(0, -2)$) {$K^{\bullet}$};
			\end{pgfonlayer}
		\end{tikzpicture}
\vspace{-3mm}
		\caption{Knots $K$, $K'$ and $K^{\bullet}$.}
		\label{fig:KK'}
	\end{figure}
	
	Then, $v(K) - v(K') = v(K^\bullet)$ according to Vassiliev skein relation. However, a knot $K^\bullet$ has $n + 1$ double points by construction, whereas $v$ is an invariant of order $\le n$, so, by Definition\,\ref{def:Vass_inv_def}, we have $v(K^\bullet) = 0$. Therefore, $v(K)=v(K')$ for all $v \in \mathcal{V}_n$ and $K$, $K'$ are singular knots with $n$ double points that differ as shown in Fig.\,\ref{fig:KK'}.
	
It follows that changing one of intersection points in a knot diagram does not affect a value of an invariant $v$ of a singular knot with $n$ double points. Hence, a value of $v \in \mathcal{V}_n$ depends only on an order in which we meet the double points when moving along the singular knot.

Information about a relative position of double points of a singular knot can be presented in form of a chord diagram, where each double point corresponds to a chord with its endpoints being preimage (with respect to an immersion that defines a singular knot) of double points of a singular knot.

A chord diagram constructed in accordance with these rules is referred to as a {\it chord diagram of a singular knot}. A chord diagram generated from a singular knot $K$ is denoted by $\sigma (K)$. Examples of this correspondence between singular knots and their chord diagrams are provided in Table\,\ref{tab:knot-chord}:
\begin{table}[!h]
		\centering
		\begin{tabular}{ r| m{1.9cm} |m{1.9cm}|m{1.9cm}|m{1.9cm}|m{1.9cm}| }
			$K$&
			\resizebox*{1.8cm}{1.8cm}{\begin{tikzpicture}
					\begin{knot}[
						style= thick,
						consider self intersections=true,
						%						draft mode=crossings,
						flip crossing=2,
						only when rendering/.style={
						}
						]
						\strand (0,2) .. controls +(2.2,0) and +(120:-2.2) .. (210:2) .. controls +(120:2.2) and +(60:2.2) .. (-30:2) .. controls +(60:-2.2) and +(-2.2,0) .. (0,2);
					\end{knot}
			\end{tikzpicture}}
			& \resizebox*{1.8cm}{1.8cm}{\begin{tikzpicture}
					
					\begin{knot}[
						style= thick,
						consider self intersections=true,
						%						draft mode=crossings,
						flip crossing=2,
						only when rendering/.style={
						}
						]
						\strand (0,2) .. controls +(2.2,0) and +(120:-2.2) .. (210:2) .. controls +(120:2.2) and +(60:2.2) .. (-30:2) .. controls +(60:-2.2) and +(-2.2,0) .. (0,2);
					\end{knot}
					\draw[fill] (149: 0.74) circle (0.14);
			\end{tikzpicture}}
			& \resizebox*{1.8cm}{1.8cm}{\begin{tikzpicture}
					
					\begin{knot}[
						style= thick,
						consider self intersections=true,
						%						draft mode=crossings,
						flip crossing=2,
						only when rendering/.style={
						}
						]
						\strand (0,2) .. controls +(2.2,0) and +(120:-2.2) .. (210:2) .. controls +(120:2.2) and +(60:2.2) .. (-30:2) .. controls +(60:-2.2) and +(-2.2,0) .. (0,2);
					\end{knot}
					\draw[fill] (31: 0.74) circle (0.14);
					\draw[fill] (149: 0.74) circle (0.14);
					
			\end{tikzpicture}}
			&\resizebox*{1.8cm}{1.8cm}{\begin{tikzpicture}
					
					\begin{knot}[
						style= thick,
						consider self intersections=true,
						%						draft mode=crossings,
						flip crossing=2,
						only when rendering/.style={
						}
						]
						\strand (0,2) .. controls +(2.2,0) and +(120:-2.2) .. (210:2) .. controls +(120:2.2) and +(60:2.2) .. (-30:2) .. controls +(60:-2.2) and +(-2.2,0) .. (0,2);
					\end{knot}
					\draw[fill] (31: 0.74) circle (0.14);
					\draw[fill] (149: 0.74) circle (0.14);
					\draw[fill] (0, -0.74) circle (0.14);
					
			\end{tikzpicture}}
			&\resizebox*{1.8cm}{1.8cm}{\begin{tikzpicture}[use Hobby shortcut]
					\begin{knot}[
						consider self intersections=true,
						%						  draft mode=crossings,
						ignore endpoint intersections=false,
						flip crossing=4,
						only when rendering/.style={
							%    show curve endpoints
						}
						]
						\strand ([closed]0,0) .. (1.5,1) .. (.5,2) .. (-.5,1) .. (.5,0) .. (0,-.5) .. (-.5,0) .. (.5,1) .. (-.5,2) .. (-1.5,1) .. (0,0);
					\end{knot}
					\path (0,-.7);
					\draw[fill] (0, 1.91) circle (0.1);
					\draw[fill] (0, .6) circle (0.1);
					\draw[fill] (0.5, .05) circle (0.1);
			\end{tikzpicture}}
			\\
			\hline & & & & &
			\\ 
			$\sigma(K)$&
			\resizebox*{1.8cm}{1.8cm}{  \begin{tikzpicture}
					\draw[ultra thick] (6,0) circle (1);
					
			\end{tikzpicture}}
			& \resizebox*{1.8cm}{1.8cm}{  \begin{tikzpicture}
					\draw[ultra thick] (6,0) circle (1);
					\draw (5, 0) -- (7,0);
			\end{tikzpicture}}
			& \resizebox*{1.8cm}{1.8cm}{  \begin{tikzpicture}
					\draw[ultra thick] (6,0) circle (1);
					\draw (5, 0) -- (7,0);
					\draw (6,-1) -- (6,1);
			\end{tikzpicture}}
			&\resizebox*{1.8cm}{1.8cm}{  \begin{tikzpicture}
					\draw[ultra thick] (0,0) circle (1);
					\draw (-1, 0) -- (1,0);
					\draw (60:1) -- (240:1);
					\draw (-60:1)--(120:1);
			\end{tikzpicture}}
			&\resizebox*{1.8cm}{1.8cm}{  \begin{tikzpicture}
					\draw[ultra thick] (0,0) circle (1);
					\draw (-1, 0) -- (1, 0);
					\draw (70:1)-- (-70: 1);
					\draw (110:1)--(250:1);
			\end{tikzpicture}}
			\\
		\end{tabular}
		\caption{Chord diagrams of singular knots.}
		\label{tab:knot-chord}
	\end{table}

	\subsection{Weight systems}
	
	In the previous paragraph it has been shown that a value of an invariant $v\in \mathcal{V}_n$ on a singular knot $K$ with $n$ double points equals to a value of some function on a chord diagram $\sigma(K)$. Formally, there exists a mapping $\alpha_n$ defined on the set $\mathcal{V}_n$, which takes values in the set of functions mapping chord diagrams with $n$ chords to some numbers.
	
	\begin{equation}
		v(K) = \alpha_n(v) (\sigma(K))\,.
	\end{equation}
	
	However, not every function on the set of chord diagrams comes from a Vassiliev invariant. To verify this, consider a singular knot and its chord diagram shown in Fig.\,\ref{fig:1T}.
	
	\begin{figure}[!h]
		\centering
		\begin{tikzpicture}[scale=0.8]
			\draw[dashed] (-4,0) rectangle (-1.5,2);
			\draw[dashed] (1.5, 0) rectangle (4, 2);
			\draw plot[smooth, tension=0.7] coordinates {(-1.5, 1.75)  (1.5, 0.25) };
			\draw (-1.5, 0.25) -- (1.5, 1.75);
			\draw[fill] (0, 1) circle (0.1);
			\node at (-2.75, 1) {tangle};
			\node at (2.75, 1) {tangle};
			\draw[dashed] (7, 1) circle (1);
			\draw[thick] (7, 1)+(160:1) arc(160:200:1);
			\draw[thick] (7, 1)+(340:1) arc(340:380:1);
			\draw (6,1)--(8,1);
			\draw[dotted] ($(7, 1)+(160:1)$) --($(7, 1)+(20:1)$);
			\draw[dotted] ($(7, 1)+(200:1)$) --($(7, 1)+(340:1)$);
		\end{tikzpicture}
		\caption{Special knot and its chord diagram.}
		\label{fig:1T}
	\end{figure}

	It follows from Vassiliev skein relation that any Vassiliev invariant of a singular knot of the type shown in Fig.\,\ref{fig:1T} must be zero. Therefore, all functions on the set of chord diagrams derived from Vassiliev invariants must vanish on a diagram with a chord that does not intersect with any other chords. This condition is known as the {\it one-term relation} \cite{chmutov2012introduction} or {\it framing-independence relation}.
	
	Another restriction can be derived by considering the expression in Fig.\,\ref{fig:4T-knot}, where $v$ is an arbitrary Vassiliev invariant and singular knots where it is evaluated coincide everywhere except for a region in which they differ as shown in Fig.\,\ref{fig:4T-knot}.

	\begin{figure}[!h]
		\centering
		\begin{tikzpicture}[scale=0.7, use Hobby shortcut]
			\coordinate (a) at (-2, 0);
			\coordinate (b) at (2, 0);
			\coordinate (c) at (6, 0);
			\coordinate (d) at (10, 0);
			\node[scale=1.6] at ($(a)-(1.6, 0)$) {$v$};
			\node[scale=2.5] at ($(a)-(1.2, 0)$) {$($};
			\node[scale=2.5] at ($(a)+(1.2, 0)$) {$)$};
			
			\node[scale=1.8] at ($(b)-(2.2, 0)$) {$-$};
			\node[scale=1.6] at ($(b)-(1.6, 0)$) {$v$};
			\node[scale=2.5] at ($(b)-(1.2, 0)$) {$($};
			\node[scale=2.5] at ($(b)+(1.2, 0)$) {$)$};
			
			\node[scale=1.8] at ($(c)-(2.2, 0)$) {$-$};
			\node[scale=1.6] at ($(c)-(1.6, 0)$) {$v$};
			\node[scale=2.5] at ($(c)-(1.2, 0)$) {$($};
			\node[scale=2.5] at ($(c)+(1.2, 0)$) {$)$};
			
			\node[scale=1.8] at ($(d)-(2.2, 0)$) {$+$};
			\node[scale=1.6] at ($(d)-(1.6, 0)$) {$v$};
			\node[scale=2.5] at ($(d)-(1.2, 0)$) {$($};
			\node[scale=2.5] at ($(d)+(1.2, 0)$) {$)$};
			\node[scale=1.8] at ($(d)+(1.9, 0)$) {$=$};
			\node[scale=1.6] at ($(d)+(2.7, 0.1)$) {$0$};
			\draw[dashed] (b) circle (1);
			
			\draw[-{Stealth[length=2.5mm]}, thick] ($(b)-(1, 0)$)--($(b)+(1,0)$);
			\draw[-{Stealth[length=2.5mm]}, thick] ($(b)+(45:0.5)$)--($(b)+(45:1)$);
			\draw[thick] ($(b)+(45:-1)$)--($(b)+(45:0.3)$);
			\draw [fill] ($(b)+(0:0.5)$) circle (0.075);
			\draw [fill] ($(b)$) circle (0.075);
			\draw[-{Stealth[length=2.5mm]}, thick] ($(b)+(90:-1)$)..($(b)+(135:-0.5)$)..($(b)+(0:0.5)$) ..($(b)+(45:0.4)$)..($(b)+(90:1)$);
			
			\draw[dashed] (a) circle (1);
			
			\draw[-{Stealth[length=2.5mm]}, thick] ($(a)-(1, 0)$)--($(a)+(1,0)$);
			\draw[-{Stealth[length=2.5mm]}, thick] ($(a)+(45:-1)$)--($(a)+(45:1)$);
			\draw [fill] ($(a)+(0:-0.5)$) circle (0.075);
			\draw [fill] ($(a)$) circle (0.075);
			\draw[thick] ($(a)+(90:-1)$)..($(a)+(85:-0.8)$)..($(a)+(55:-0.5)$);
			\draw[-{Stealth[length=2.5mm]}, thick] ($(a)+(35:-0.5)$)..($(a)+(0:-0.5)$)..($(a)+(135:0.5)$) ..($(a)+(90:1)$);

			\draw[dashed] (c) circle (1);
			
			\draw[-{Stealth[length=2.5mm]}, thick] ($(c)-(1, 0)$)--($(c)+(1,0)$);
			\draw[-{Stealth[length=2.5mm]}, thick] ($(c)+(45:-1)$)--($(c)+(45:1)$);
			\draw [fill] ($(c)+(45:0.4)$) circle (0.075);
			\draw [fill] ($(c)$) circle (0.075);
			\draw[thick] ($(c)+(90:-1)$)..($(c)+(135:-0.4)$)..($(c)+(170:-0.4)$);
			\draw[-{Stealth[length=2.5mm]}, thick] ($(c)+(10:0.4)$)..($(c)+(45:0.4)$)..($(c)+(85:0.75)$)..($(c)+(90:1)$);

			\draw[dashed] (d) circle (1);
			
			\draw[-{Stealth[length=2.5mm]}, thick] ($(d)-(0.3, 0)$)--($(d)+(1,0)$);
			\draw[ thick] ($(d)-(1, 0)$)--($(d)+(-0.5,0)$);
			\draw[-{Stealth[length=2.5mm]}, thick] ($(d)+(45:-1)$)--($(d)+(45:1)$);
			\draw [fill] ($(d)+(45:-0.4)$) circle (0.075);
			\draw [fill] (d) circle (0.075);
			\draw[-{Stealth[length=2.5mm]}, thick] ($(d)+(90:-1)$)..($(d)+(45:-0.4)$)..($(d)+(0:-0.4)$)..($(d)+(100:0.8)$)..($(d)+(90:1)$);
			
		\end{tikzpicture}
		\caption{Four terms knot relation.}
		\label{fig:4T-knot}
	\end{figure}
	
	If we resolve double points according to Vassiliev skein relation, then all terms are reduced. This results in the {\it four-term ratio} in Fig.\,\ref{fig:4T-knot}. From this relation on Vassiliev invariants, there follows a similar four-term relation on the functions $f\in\text{Im}\, \alpha_n$ originating from Vassiliev invariants (Fig.\,\ref{fig:4T-fun}) \cite[\S 4.2]{chmutov2012introduction}:

	\begin{figure}[!h]
		\centering
		%		\resizebox*{\linewidth}{2cm}
		{
			\begin{tikzpicture}[scale=0.7]
				\coordinate (a) at (-2, 0);
				\coordinate (b) at (2, 0);
				\coordinate (c) at (6, 0);
				\coordinate (d) at (10, 0);
				\node[scale=1.6] at ($(a)-(1.6, 0)$) {$f$};
				\node[scale=2.5] at ($(a)-(1.2, 0)$) {$($};
				\node[scale=2.5] at ($(a)+(1.2, 0)$) {$)$};
				
				\node[scale=1.8] at ($(b)-(2.2, 0)$) {$-$};
				\node[scale=1.6] at ($(b)-(1.6, 0)$) {$f$};
				\node[scale=2.5] at ($(b)-(1.2, 0)$) {$($};
				\node[scale=2.5] at ($(b)+(1.2, 0)$) {$)$};
				
				\node[scale=1.8] at ($(c)-(2.2, 0)$) {$-$};
				\node[scale=1.6] at ($(c)-(1.6, 0)$) {$f$};
				\node[scale=2.5] at ($(c)-(1.2, 0)$) {$($};
				\node[scale=2.5] at ($(c)+(1.2, 0)$) {$)$};
				
				\node[scale=1.8] at ($(d)-(2.2, 0)$) {$+$};
				\node[scale=1.6] at ($(d)-(1.6, 0)$) {$f$};
				\node[scale=2.5] at ($(d)-(1.2, 0)$) {$($};
				\node[scale=2.5] at ($(d)+(1.2, 0)$) {$)$};
				\node[scale=1.8] at ($(d)+(1.9, 0)$) {$=$};
				\node[scale=1.6] at ($(d)+(2.6, 0.1)$) {$0$};
				\draw[dashed] (a) circle (1);
				\draw[thick] ($(a)+(240:1)$) arc(240:300:1);
				\draw[thick] ($(a)+(0:1)$) arc(0:60:1);
				\draw[thick] ($(a)+(140:1)$) arc(140:180:1);
				\draw ($(a)+(255:1)$)--($(a)+(160:1)$);
				\draw ($(a)+(280:1)$)--($(a)+(20:1)$);
				\draw[dashed] (b) circle (1);
				\draw[thick] ($(b)+(240:1)$) arc(240:300:1);
				\draw[thick] ($(b)+(0:1)$) arc(0:60:1);
				\draw[thick] ($(b)+(140:1)$) arc(140:180:1);
				\draw ($(b)+(285:1)$)--($(b)+(160:1)$);
				\draw ($(b)+(255:1)$)--($(b)+(20:1)$);
				\draw[dashed] (c) circle (1);
				\draw[thick] ($(c)+(240:1)$) arc(240:300:1);
				\draw[thick] ($(c)+(0:1)$) arc(0:60:1);
				\draw[thick] ($(c)+(140:1)$) arc(140:180:1);
				\draw ($(c)+(45:1)$)--($(c)+(160:1)$);
				\draw ($(c)+(270:1)$)--($(c)+(20:1)$);
				\draw[dashed] (d) circle (1);
				\draw[thick] ($(d)+(240:1)$) arc(240:300:1);
				\draw[thick] ($(d)+(0:1)$) arc(0:60:1);
				\draw[thick] ($(d)+(140:1)$) arc(140:180:1);
				\draw ($(d)+(20:1)$)--($(d)+(160:1)$);
				\draw ($(d)+(270:1)$)--($(d)+(45:1)$);
		\end{tikzpicture}}
		\caption{4-term relation on functions coming from Vassiliev invariants.}
		\label{fig:4T-fun}
	\end{figure}
	
	\noindent The derived conditions, namely, the one-term and the four-term relations, serve as a basis for the following definitions.
	
	\begin{definition}[\textbf{Space of chord diagrams $\mathcal{A}$} \cite{bar1995vassiliev}]
		The space of chord diagrams $\mathcal{A}_n$ is a linear space generated by all chord diagrams with $n$ chords factorized by a subspace generated by four-term elements of the form:
		\begin{figure*}[!h]
			\centering
			\begin{tikzpicture}[scale=0.5]
				\coordinate (a) at (-2, 0);
				\coordinate (b) at (1, 0);
				\coordinate (c) at (4, 0);
				\coordinate (d) at (7, 0);
				\draw[dashed] (a) circle (1);
				\draw[thick] ($(a)+(240:1)$) arc(240:300:1);
				\draw[thick] ($(a)+(0:1)$) arc(0:60:1);
				\draw[thick] ($(a)+(140:1)$) arc(140:180:1);
				\draw ($(a)+(255:1)$)--($(a)+(160:1)$);
				\draw ($(a)+(280:1)$)--($(a)+(20:1)$);
				
				\draw[dashed] (b) circle (1);
				\draw[thick] ($(b)+(240:1)$) arc(240:300:1);
				\draw[thick] ($(b)+(0:1)$) arc(0:60:1);
				\draw[thick] ($(b)+(140:1)$) arc(140:180:1);
				\draw ($(b)+(285:1)$)--($(b)+(160:1)$);
				\draw ($(b)+(255:1)$)--($(b)+(20:1)$);
				
				\draw[dashed] (c) circle (1);
				\draw[thick] ($(c)+(240:1)$) arc(240:300:1);
				\draw[thick] ($(c)+(0:1)$) arc(0:60:1);
				\draw[thick] ($(c)+(140:1)$) arc(140:180:1);
				\draw ($(c)+(45:1)$)--($(c)+(160:1)$);
				\draw ($(c)+(270:1)$)--($(c)+(20:1)$);
				
				\draw[dashed] (d) circle (1);
				\draw[thick] ($(d)+(240:1)$) arc(240:300:1);
				\draw[thick] ($(d)+(0:1)$) arc(0:60:1);
				\draw[thick] ($(d)+(140:1)$) arc(140:180:1);
				\draw ($(d)+(20:1)$)--($(d)+(160:1)$);
				\draw ($(d)+(270:1)$)--($(d)+(45:1)$);
				
				\node at ($0.5*(a)+0.5*(b)$) {$-$};
				\node at ($0.5*(c)+0.5*(b)$) {$-$};
				\node at ($0.5*(c)+0.5*(d)$) {$+$};
			\end{tikzpicture}
		\end{figure*}
		\newline
		and by one-term elements of the form:
		
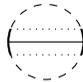
\begin{figure*}[!h]
			\centering
			\begin{tikzpicture}[scale=0.5]
				\draw[dashed] (7, 1) circle (1);
				\draw[thick] (7, 1)+(160:1) arc(160:200:1);
				\draw[thick] (7, 1)+(340:1) arc(340:380:1);
				\draw (6,1)--(8,1);
				\draw[dotted] ($(7, 1)+(160:1)$) --($(7, 1)+(20:1)$);
				\draw[dotted] ($(7, 1)+(200:1)$) --($(7, 1)+(340:1)$);
			\end{tikzpicture}
		\end{figure*}\\
		
		\noindent The whole space of chord diagrams ${\cal A}=\bigoplus\limits_n \mathcal{A}_n\,$.
	\end{definition}
	
	More compactly, we can write that $\mathcal{A}_n = \textbf{A}_n/\langle\text{4T, 1T} \rangle$, where 4T and 1T are the four-term and the one-term relations correspondingly. 
	
	\begin{definition}[\textbf{Weight systems} $\mathcal{W}$ \cite{bar1995vassiliev}]
		The space of weight systems $\mathcal{W}_n = \mathcal{A}_n^{*} = \text{Hom}~(\mathcal{A}_n, \mathbb{R})$ is a space of linear functions on $\mathcal{A}_n$. 
		\label{def:W_def}
	\end{definition}

	Due to these definitions, the weight systems satisfy the one-term and four-term relations, which are conditions for functions derived from Vassiliev invariants. 
	
	\subsection{Vassiliev-Kontsevich theorem}\label{sec:VK-theorem}
	In the previous paragraph, we have shown that any Vassiliev invariant generates a weight system, which is a linear mapping $\alpha_n:\mathcal{V}_n \rightarrow \mathcal{W}_n$. Due to the definition of Vassiliev invariants, the kernel of this mapping is $\text{Ker}(\alpha_n) = \mathcal{V}_{n - 1}\,$. This implies that there exists an injective operator $\bar\alpha_n : \mathcal{V}_n / \mathcal{V}_{n - 1} \rightarrow \mathcal{W}_n$ (see \cite{chmutov2012introduction}).

	M. Kontsevich proved in \cite{kontsevich1993vassiliev} that $\bar{\alpha}_n$ is an isomorphism. In other words, any weight system $w \in \mathcal{W}_n$ can be generated by a Vassiliev invariant $v \in \mathcal{V}_n / \mathcal{V}_{n - 1}$. A value of this invariant on a knot $K$ is given by the following formula:
	
	\begin{equation}
		v(K) = w(I(K))\,.
	\end{equation}
	Here, $I(K)$ is the \textit{Kontsevich integral} of 	the knot $K$. Its formal definition can be found in \cite[Chapter 8]{chmutov2012introduction}. A value of the Kontsevich integral of a knot $K$ is a formal series of chord diagrams $D_{n,m}$ in the space $\mathcal{A}_n$, with the coefficients of this series being finite-order invariants $v_{n,m}(K)$:
	\be
	\label{ki}
	I(K) = \sum_{n=0}^{\infty} \hbar^n \sum_{m=1}^{\text{dim} \mathcal{A}_n} v_{n,m}(K)\, D_{n,m}\,,
	\ee
	where $\hbar$ is a formal variable. Furthermore, all Vassiliev invariants can be derived from the coefficients of the Kontsevich integral, that is why it is also referred to as the {\it universal Vassiliev invariant} \cite{chmutov2012introduction}.

	Thus, \textit{Vassiliev-Kontsevich theorem} \cite[\S 4.2.]{chmutov2012introduction} is valid:
	
	\begin{equation}
		\mathcal{W} = \bigoplus_{n=0}^{\infty} \mathcal{W}_n \cong \bigoplus_{n=0}^{\infty} \mathcal{V}_n/\mathcal{V}_{n-1} = \mathcal{V}\,.
	\end{equation}

	\subsection{Jacobi diagrams}\label{sec:Jacobi-diag}
	
	\begin{definition}[\textbf{Space of Jacobi diagrams} $\mathcal{C}$ \cite{bar1995vassiliev, chmutov2012introduction}]
		The space of (closed) Jacobi diagrams $\mathcal{C}_n$ is a linear space generated by connected graphs having $2n$ trivalent vertices, as well as having a distinguished cycle (vertices on this cycle are called external and the cycle is called \textit{Wilson loop}), the remaining (internal) vertices are equipped with a cyclic order of half-edges, modulo the STU relation:
		\begin{figure}[!h]
			\centering
			\begin{tikzpicture}[scale=0.7]
				\coordinate (s) at (0,0);
				\draw[-{Stealth[length=3mm]}, ultra thick] ($(s)+ (45:-2)$) arc(225:315:2);
				\coordinate (t) at (4,0);
				\draw[-{Stealth[length=3mm]}, ultra thick] ($(t)+ (45:-2)$) arc(225:315:2);
				\coordinate (u) at (8,0);
				\draw[-{Stealth[length=3mm]}, ultra thick] ($(u)+ (45:-2)$) arc(225:315:2);
				\coordinate (s1) at ($(s)+(-1, 0)$);
				\coordinate (t1) at ($(t)+(-1, 0)$);
				\coordinate (u1) at ($(u)+(-1, 0)$);
				\coordinate (s2) at ($(s)+(1, 0)$);
				\coordinate (t2) at ($(t)+(1, 0)$);
				\coordinate (u2) at ($(u)+(1, 0)$);
				\draw ($(s)+(0, -2)$) -- ($(s)+(0, -1)$);
				\draw (s1) -- ($(s)+(0, -1)$);
				\draw ($(s)+(0, -1)$) -- (s2);
				\draw[fill] ($(s)+(0, -1)$) circle(0.1);
				
				\draw (t1) -- ($(t)+(70:-2)$);
				\draw ($(t)+(110:-2)$) -- (t2);
				
				\draw (u1) -- ($(u)+(110:-2)$);
				\draw ($(u)+(70:-2)$) -- (u2);
				\node[scale=2.5] at ($0.5*(s) + 0.5*(t)-(0,1)$) {$=$};
				\node[scale=2.5] at ($0.5*(u) + 0.5*(t)-(0, 1)$) {$-$};
			\end{tikzpicture}
			\caption{STU relation.}
			\label{fig:STU-relation}
		\end{figure}\\

		\noindent The whole space of Jacobi diagrams ${\cal C}= \bigoplus\limits_n \mathcal{C}_n$.
	\end{definition}

	 Examples of Jacobi diagrams are presented in Fig.\,\ref{fig:Jacobi_examples}. In this and subsequent discussions, we assume that the cyclic ordering of the half-edges at each vertex is consistent with the direction of the drawing surface. Special consideration should be given to diagram \textbf{d} in Fig.\ref{fig:Jacobi_examples}, which is equal to zero as an element of $\mathcal{C}$ due to the STU relation.
	
	\begin{figure}[!h]
		\centering
		\begin{tikzpicture}[scale=0.8]
			\draw[ultra thick]  (-2, 0) circle (1);
			\draw (-2, 0) -- ($(-2, 0)+(0:1)$);
			\draw (-2, 0) -- ($(-2, 0)+(120:1)$);
			\draw (-2, 0) -- ($(-2, 0)+(240:1)$);
			\draw[fill] (-2, 0) circle  (0.05);
			\draw[ultra thick]  (2, 0) circle (1);
			\draw ($(2, 0)+(45:0.5)$) -- ($(2, 0)+(45:1)$);
			\draw ($(2, 0)+(135:0.5)$) -- ($(2, 0)+(135:1)$);
			\draw ($(2, 0)+(225:0.5)$) -- ($(2, 0)+(225:1)$);
			\draw ($(2, 0)+(315:0.5)$) -- ($(2, 0)+(315:1)$);
			\draw (2, 0) circle (0.5);
			\draw[fill] ($(2, 0)+(315:0.5)$) circle  (0.05);
			\draw[fill] ($(2, 0)+(225:0.5)$) circle  (0.05);
			\draw[fill] ($(2, 0)+(135:0.5)$) circle  (0.05);
			\draw[fill] ($(2, 0)+(45:0.5)$) circle  (0.05);
			\draw[ultra thick] (6,0) circle (1);
			\draw (5, 0) -- (7,0);
			\draw (6,-1) -- (6,1);
			\draw[ultra thick]  (10, 0) circle (1);
			\draw (9, 0) -- ($(10, 0)+(180:0.4)$);
			\draw (10, 0) circle (0.4);
			\draw[fill] (9.6, 0) circle  (0.05);
			
			\node at (-2, -1.5) {\textbf{a}};
			
			\node at (2, -1.5) {\textbf{b}};
			
			\node at (6, -1.5) {\textbf{c}};
			
			\node at (10, -1.5) {\textbf{d}};
		\end{tikzpicture}
		\caption{Examples of Jacobi diagrams.}
		\label{fig:Jacobi_examples}
	\end{figure}
	
	Chord diagrams are a special type of Jacobi diagrams. In addition, the four-term relation follows from the STU relation, so it automatically holds in $\mathcal{C}$. On the other hand, by consistently applying the STU relation, any Jacobi diagram can be reduced to a linear combination of chord diagrams. This linear combination generally depends on an order in which the STU relations are applied, but all such combinations are equivalent modulo the 4T relation. The statement about this fact can be found in \cite[\S 5.3.]{chmutov2012introduction}.

	Thus, $\textbf{A}_n/  \langle \text{4T}\rangle \cong \mathcal{C}_n$, and $\mathcal{A}_n \cong \mathcal{C}_n / \langle \text{1T}\rangle$. This fact allows us to work with Jacobi diagrams instead of chord diagrams. The advantage of Jacobi diagrams is that the STU relation is more convenient than the four-term relation.
	
	In the space of Jacobi diagrams, there is a multiplication operation. The product of two Jacobi diagrams is a new diagram that is created by cutting  the Wilson loops from the multiplier diagrams and then gluing them together. An example of this product is shown in Fig. \ref{fig:diag-mult}.
	
	\begin{figure}[!h]
		\centering
		\begin{tikzpicture}[scale=0.8]
			\coordinate (a) at (-3, 0);
			\coordinate (b) at (0, 0);
			\coordinate (c) at (3, 0);
			\draw[ultra thick] (a) circle (1);
			\draw[ultra thick] (b) circle (1);
			\draw[ultra thick] (c) circle (1);
			\node [scale=2] at ($0.5*(a)+0.5*(b)$) {$\cdot$};
			\node at ($0.5*(b)+0.5*(c)$) {$=$};
			
			\draw (a) -- +(60:1);
			\draw (a) -- + (180:1);
			\draw (a) -- +(-60:1);
			\draw[fill] (a) circle (0.05);
			
			\draw (b) circle (0.4);
			\draw (b) +(45:0.4)-- +(45:1);
			\draw (b) +(135:0.4)-- +(135:1);
			\draw (b) +(225:0.4)-- +(225:1);
			\draw (b) +(315:0.4)-- +(315:1);
			\draw[fill] (b)+(45:0.4) circle (0.05);
			\draw[fill] (b)+(135:0.4) circle (0.05);
			\draw[fill] (b)+(225:0.4) circle (0.05);
			\draw[fill] (b)+(315:0.4) circle (0.05);
			
			\draw (c)+(0,0.3) circle (0.3);
			\draw (c)+($(0,0.3)+(45:0.3)$)-- +(60:1);
			\draw (c)+($(0,0.3)+(135:0.3)$)--+(120:1);
			\draw (c)+($(0,0.3)+(225:0.3)$)--+(180:1);
			\draw (c)+($(0,0.3)+(315:0.3)$)--+(0:1);
			\draw[fill] (c)+($(0,0.3)+(45:0.3) $)circle (0.05);
			\draw[fill] (c)+($(0,0.3)+(135:0.3)$) circle (0.05);
			\draw[fill] (c)+($(0,0.3)+(225:0.3)$) circle (0.05);
			\draw[fill] (c)+($(0,0.3)+(315:0.3)$) circle (0.05);
			
			\draw[fill] (c) + (0, -0.5) circle (0.05);
			
			\draw ($(c)+(-30:1)$)--($(c)+(30:-1)$);
			\draw (c) +(0, -0.5)-- +(0,-1);
			
		\end{tikzpicture}
		\caption{Product of Jacobi diagrams.}
		\label{fig:diag-mult}
	\end{figure}
	
	The space of \textit{primitive} Jacobi diagrams $\mathcal{P}$ is a subspace of $\mathcal{C}$ linearly generated by Jacobi diagrams with a connected internal graph. More strictly, the subspace of primitive diagrams is defined with respect to the co-multiplication operation in $\mathcal{C}$ (see \cite[\S 5.5.]{chmutov2012introduction}). Any Jacobi diagram can be represented as a polynomial in primitive diagrams.

\subsubsection{AS and IHX relations}

The STU relation imposes additional constraints on an internal graph of Jacobi diagrams. This can be seen in the example of diagram {\bf d} in Fig.\,\ref{fig:Jacobi_examples}. It can be shown that all Jacobi diagrams having such a loop vanish by STU relations.

This statement is a special case of the AS relation shown in Fig.\,\ref{fig:AS-relation}. 

\begin{figure}[!h]
	\centering
	\begin{tikzpicture}[use Hobby shortcut, scale=0.7]
		\coordinate (a) at (-2, 0);
		\coordinate (b) at (2, 0);
		\draw [fill] (a) circle (0.05);
		\draw [fill] (b) circle (0.05);
		\draw (a) -- + (90: -1);
		\draw (a) .. +(45: 0.2) ..+(90:0.4) ..+(120:0.7) ..+(125: 1.2)..+(120: 2);
		\draw (a) .. +(135: 0.2) .. +(90:0.4) ..+(60:0.7) ..+(55: 1.2)..+(60: 2);
		
		\node [scale=1.5] at (0,0.2) {$= \;-$};
		\draw (b) -- +(90: -1);
		\draw (b) -- +(60: 2);
		\draw (b) -- +(120: 2);
	\end{tikzpicture}
	\caption{AS relation.}
	\label{fig:AS-relation}
\end{figure}

There is also the second relation called the IHX relation, see Fig.\,\ref{fig:IHX-relation}.

\begin{figure}[!h]
	\centering
	\begin{tikzpicture}[scale=0.9]
		\coordinate (d) at (3., 0);
		\coordinate (c) at (2., 0);
		\draw (c)--(d);
		\draw [fill] (c) circle (0.05);
		\draw [fill] (d) circle (0.05);
		\draw (c) -- ($(c)+(120:1)$);
		\draw (c) -- ($(c)+(240:1)$);	
		\draw (d) -- ($(d)+(120:-1)$);
		\draw (d) -- ($(d)+(240:-1)$);

		\coordinate (e) at (-1, 0.5);
		\coordinate (f) at (-1, -0.5);
		\draw [fill] (f) circle (0.05);
		\draw [fill] (e) circle (0.05);
		\draw (e)--(f);
		\draw (e) -- +(30:1);
		\draw (e) -- +(150:1);
		\draw (f) -- +(150:-1);
		\draw (f) -- +(30:-1);
		\node [scale=1.5] at ($0.2*(e)+0.2*(f)+0.6*(c)$) {$=$};

		\coordinate (a) at (6., 0);
		\coordinate (b) at (5., 0);
		\draw (a)--(b);
		\draw [fill] (a) circle (0.05);
		\draw [fill] (b) circle (0.05);
		\draw (a) -- ($(a)+(-60:1)$);
		\draw (a) -- ($(a)+(120:1)+(-1, 0)$);	
		\draw (b) -- ($(b)+(60:-1)$);
		\draw (b) -- ($(b)+(60:1)+(1, 0)$);	
		\node at ($0.5*(b) + 0.5*(d)$) {$-$};
	\end{tikzpicture}

	\caption{IHX relation.}
	\label{fig:IHX-relation}
\end{figure}

The relations AS and IHX hold true with respect to the internal graphs of diagrams in the space $\mathcal{C}$. These relations can be derived from the STU relation. For more details on this statement, see \cite[\S 5.2]{chmutov2012introduction}.

\subsubsection{Open Jacobi diagrams}\label{sec:open_Jacobi}

	\begin{definition}[Open Jacobi diagrams $\mathcal{B}$ \cite{bar1995vassiliev}]
	 The space of open Jacobi diagrams $\mathcal{B}$ is generated by trivalent diagrams modulo the AS and IHX relations with unordered (unfixed) legs and with at least one leg in each connected component. The legs are also referred to as univalent vertices.
	 
	 \noindent The space $\mathcal{B}$ admits two gradings: by the number of legs $l$ and by half the total number of vertices (both trivalent and univalent) $n$. 
	 
\end{definition}

The map $\rho$ from $\mathcal{B}$ to $\mathcal{C}$ is defined as shown in Fig.\,\ref{fig:b-c-morph}. Open Jacobi diagrams are {\it unfixed diagrams} as they do not carry any order of their legs. This is why the map $\rho$ averages between all the different ways of attaching legs to the Wilson loop.

\begin{figure}[!h]
	\centering
	\begin{tikzpicture}[scale=1.5]
		\coordinate (t) at (0, 0);
		\coordinate (a) at ($(t)+(0, 0.2)$);
		\node[scale=1.4] at ($(t)+(-0.75, 0)$) {$\rho$};
		\node[scale=2.5] at ($(t)+(-0.52, 0)$) {$($};
		\node[scale=2.5] at ($(t)+(+0.52, 0)$) {$ ) $};
		\draw[dotted] (a) circle (0.3);
		\draw (a) +(50:-0.3)--+(-0.2, -0.6);
		\draw (a) +(70:-0.3)--+(-0.1, -0.6);
		\draw (a) +(90:-0.3)--+(0, -0.6);
		\draw (a) +(110:-0.3)--+(0.1, -0.6);
		\draw (a) +(130:-0.3)--+(0.2, -0.6);
		
		\coordinate (b) at (3.3, 0);

		\node[scale=1.5] at ($0.5*(b)+0.5*(t)+(-0.1, -0.14)$) {$=\frac{1}{n!} \sum\limits_{s \in S_n}$};
		%			\node[scale=1.5] at ($(b)+(0.82, 0)$) {$)$};
		\draw[thick] (b) circle (0.7);
		\draw[dotted] (b)+(0, 0.3) circle (0.3);
		\draw[thin] (b)+(-0.3, -0.5) rectangle +(0.3, -0.15);
		\node[scale=0.8] at ($(b)+(0, -0.33)$) {$s$};
		\draw ($(b)+(0, 0.3)+(50:-0.3)$) -- ($(b)+(-0.2,-0.15)$);
		\draw ($(b)+(0, 0.3)+(70:-0.3)$) -- ($(b)+(-0.1,-0.15)$);
		\draw ($(b)+(0, 0.3)+(90:-0.3)$) -- ($(b)+(0,-0.15)$);
		\draw ($(b)+(0, 0.3)+(110:-0.3)$) -- ($(b)+(0.1,-0.15)$);
		\draw ($(b)+(0, 0.3)+(130:-0.3)$) -- ($(b)+(0.2,-0.15)$);
		
		\draw ($(b)+(70:-0.7)$) -- ($(b)+(-0.2,-0.5)$);
		\draw ($(b)+(80:-0.7)$) -- ($(b)+(-0.1,-0.5)$);
		\draw ($(b)+(90:-0.7)$) -- ($(b)+(0,-0.5)$);
		\draw ($(b)+(100:-0.7)$) -- ($(b)+(0.1,-0.5)$);
		\draw ($(b)+(110:-0.7)$) -- ($(b)+(0.2,-0.5)$);
	\end{tikzpicture}
	\caption{Isomorphism $\rho: \mathcal{B} \rightarrow \mathcal{C}$.}
	\label{fig:b-c-morph}
\end{figure}

The proof of $\rho: \mathcal{B} \rightarrow \mathcal{C}$ being isomorphism, can be found in \cite[\S 5.7]{chmutov2012introduction}. The benefit of operating in $\mathcal{B}$ is an extra grading by the number of legs, which is very convenient in our study in Section \ref{sec:weight-system-kernel}

\subsection{Lie algebra weight systems}\label{sec:Lie-weights}
Lie algebras give rise to linear functions on $\mathcal{C}$ (weight systems). For a simple Lie algebra $L$ with a nondegenerate Killing form and its representation $R$ the construction is the following. External vertices of a Jacobi diagram are mapped to generators $T_R^a $ of an algebra $L$ in a representation $R$. Internal vertices correspond to structure constants $f^{a b c}$, internal edges correspond to the  metric tensor $ g_{ab}$. The Wilson loop corresponds to taking the trace of the product of generators divided by the dimension of the representation $R$.

The above construction belongs to Bar-Natan \cite{bar1991weights, bar1995vassiliev} and is referred to as the weight system of a Lie algebra. The consistency with the STU relation is ensured by the following equation: $T^a T^b - T^b T^a = [T^a, T^b] = {f^{ab}}_c T^c$. An illustration of this equation can be seen in Fig.\,\ref{fig:STU-Lie}.

\begin{figure}[!h]
	\centering
	\begin{tikzpicture}[scale=0.7]
		\coordinate (s) at (0,0);
		\draw[-{Stealth[length=3mm]}, ultra thick] ($(s)+ (45:-2)$) arc(225:315:2);
		\coordinate (t) at (4,0);
		\draw[-{Stealth[length=3mm]}, ultra thick] ($(t)+ (45:-2)$) arc(225:315:2);
		\coordinate (u) at (8,0);
		\draw[-{Stealth[length=3mm]}, ultra thick] ($(u)+ (45:-2)$) arc(225:315:2);
		\node (s1) at ($(s)+(-1, 0)$) {$a$};
		\node (t1) at ($(t)+(-1, 0)$) {$a$};
		\node (u1) at ($(u)+(-1, 0)$) {$a$};
		\node (s2) at ($(s)+(1, 0)$) {$b$};
		\node (t2) at ($(t)+(1, 0)$) {$b$};
		\node (u2) at ($(u)+(1, 0)$) {$b$};
		\draw ($(s)+(0, -2)$) -- ($(s)+(0, -1)$);
		\draw (s1) -- ($(s)+(0, -1)$);
		\draw ($(s)+(0, -1)$) -- (s2);
		\draw[fill] ($(s)+(0, -1)$) circle(0.1);
		
		\draw (t1) -- ($(t)+(70:-2)$);
		\draw ($(t)+(110:-2)$) -- (t2);
		
		\draw (u1) -- ($(u)+(110:-2)$);
		\draw ($(u)+(70:-2)$) -- (u2);
		\node[scale=1.5] at ($0.5*(s) + 0.5*(t)-(0,1)$) {$=$};
		\node[scale=1.5] at ($0.5*(u) + 0.5*(t)-(0, 1)$) {$-$};

		\node at ($(s)+ (0, -2.5)$) {${f^{ab}}_c T^c$};
		
		\node at ($(t)+(70:-2)-(0, 0.5)$) {$T^a$};
		\node at ($(t)+(110:-2)-(0, 0.5)$) {$T^{b}$};
		
		\node at ($(u)+(70:-2)-(0, 0.5)$) {$T^b$};
		\node at ($(u)+(110:-2)-(0, 0.5)$) {$T^{a}$};
	\end{tikzpicture}
	\caption{STU relation as commutation in Lie algebra representation.}
	\label{fig:STU-Lie}
\end{figure} 

The weight system for an algebra $L$ and a representation $R$ is denoted\footnote{In what follows, we do not specify the representation $R$, i.e. work with an arbitrary representation, and omit the superscript $R$ in the denotation of a Lie algebra weight system $\varphi_L^R$.} $\varphi_L^R$. The factor $1/\dim R $ in the definition is necessary to ensure that a value of this weight system on a diagram without edges is $1$ (in this case, the trace is evaluated from the unit matrix with dimensions ${\dim R}\times{\dim R}$). The weight system is then multiplicative, that is, $\varphi_L^R(D_1 \cdot D_2) = \varphi^R_L(D_1)\varphi^R_L(D_2)$, where $D_1$ and $D_2$ are Jacobi diagrams.

The AS and IHX relations acquire a simple interpretation in terms of the weight systems of Lie algebras. Specifically, the AS relation is responsible for the anti-symmetry of the Lie bracket, $[x, y] = -[y, x]$ $\forall x,y\in L$, which, in the language of structure constants, can be expressed as ${f^{ab}}_c = -{f^{ba}}_c$. The IHX relation is equivalent to the Jacobi identity, $[[x, y], z] + [[y, z], x] + [[z, x], y] = 0$ $\forall x,y,z\in L$, or ${f^{ad}}_e {f^{bc}}_d + {f^{bd}}_e {f^{ca}}_d + {f^{cd}}_e {f^{ab}}_d = 0$ in terms of structure constants (see \cite[\S  5.2]{chmutov2012introduction}).

\setcounter{equation}{0}
\section{Vogel algebra}\label{sec:Vogel-algebra}
In this section, we describe main features of the framework developed by Pierre Vogel in his paper \cite{vogel2011algebraic}. Advanced features of this framework are given in Appendix~\ref{appendix}.

\begin{definition}[{\bf Algebra $\Lambda$} \cite{vogel2011algebraic}]
	$\Lambda$ is an algebra over $\mathbb{Q}$ generated by 3-legged fixed (i.e. the legs are numbered) diagrams modulo AS and IHX relations antisymmetric with respect to permutations of legs. Multiplication in $\Lambda$ is given by insertion of one diagram into any vertex of the other diagram.
	% Give definition of a fixed diagram
	\noindent $\Lambda$ is a graded algebra. Its grading is given by half the number of its trivalent vertices minus one half.
\end{definition}

\begin{remark}
	In what follows, all of the diagrams are fixed unless they are open Jacobi diagrams
\end{remark}

The condition of antisymmetry with respect to permutations of three legs in the definition of $\Lambda$ allows for insertion of a diagram $\hat{v}\in \Lambda$ into a 3-valent vertex with AS-relation. To prove that the multiplication is well-defined, it is needed to show that the product is independent of a vertex to insert $\hat{v}$ into. To do so, one needs to prove the equality between two  diagrams depicted in Fig.\,\ref{fig:2diagrams}. 

\begin{figure}[!h]
	\centering 
	\begin{tikzpicture}[scale=0.75]
		\node (a) at (-3, 0) {$\hat{v}$};
		\coordinate (b) at (-1, 0);
		\draw [fill] (b) circle (0.1);
		\draw (a) -- (b);		
		\draw (a) -- ($(a)+(120:1)$);
		\draw (a) -- ($(a)+(240:1)$);	
		\draw (b) -- ($(b)+(120:-1)$);
		\draw (b) -- ($(b)+(240:-1)$);	
		\node (d) at (4.5, 0) {$\hat{v}$};
		\coordinate (c) at (2.5, 0);
		\draw (c)--(d);
		\draw [fill] (c) circle (0.1);
		\draw (c) -- ($(c)+(120:1)$);
		\draw (c) -- ($(c)+(240:1)$);	
		\draw (d) -- ($(d)+(120:-1)$);
		\draw (d) -- ($(d)+(240:-1)$);	
		\node[scale=1.5] at ($0.5*(b) + 0.5*(c)$) {$=$};
	\end{tikzpicture}
	\caption{Two equivalent diagrams. } 
	\label{fig:2diagrams}
\end{figure}

To prove this relation the generalized form of the IHX relation is needed \cite[Lemma 5.2.9.]{chmutov2012introduction}. Applying the generalized IHX relation we end up with:

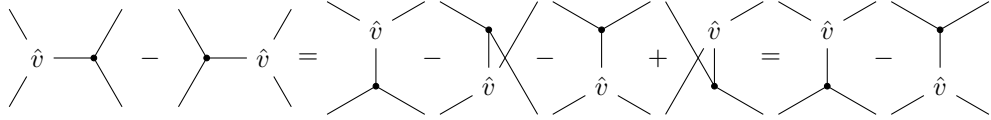
\begin{figure*}[!h]
	\centering
	\begin{tikzpicture}[scale=0.75]
		\node (a) at (-3, 0) {$\hat{v}$};
		\coordinate (b) at (-2, 0);
		\draw [fill] (b) circle (0.05);
		\draw (a) -- (b);		
		\draw (a) -- ($(a)+(120:1)$);
		\draw (a) -- ($(a)+(240:1)$);	
		\draw (b) -- ($(b)+(120:-1)$);
		\draw (b) -- ($(b)+(240:-1)$);	
		\node (d) at (1., 0) {$\hat{v}$};
		\coordinate (c) at (0., 0);
		\draw (c)--(d);
		\draw [fill] (c) circle (0.05);
		\draw (c) -- ($(c)+(120:1)$);
		\draw (c) -- ($(c)+(240:1)$);	
		\draw (d) -- ($(d)+(120:-1)$);
		\draw (d) -- ($(d)+(240:-1)$);	
		\node[] at ($0.5*(b) + 0.5*(c)$) {$-$};
		
		\node (e) at (3, 0.5) {$\hat{v}$};
		\coordinate (f) at (3, -0.5);
		\draw [fill] (f) circle (0.05);
		\draw (e)--(f);
		\draw (e) -- +(30:1);
		\draw (e) -- +(150:1);
		\draw (f) -- +(150:-1);
		\draw (f) -- +(30:-1);
		\node at ($0.2*(e)+0.2*(f)+0.6*(d)$) {$=$};
		
		\coordinate (g) at (5, 0.5);
		\node (h) at (5, -0.5) {$\hat{v}$};
		\draw (g)--(h)--+(30:-1);
		\draw (g) -- +(150:1);
		\draw (g) -- ($(g)+(0, -1)+(150:-1)$);
		\draw (h) -- ($(h)+(0, 1)+(30:1)$);
		\draw[fill] (g) circle (0.05);
		\node[] at ($0.5*(e) + 0.5*(h)$) {$-$};

		\node (l) at (7, -0.5) {$\hat{v}$};
		\coordinate (k) at (7, +0.5);
		\draw [fill] (k) circle (0.05);
		\draw (k)--(l);
		\draw (k) -- +(30:1);
		\draw (k) -- +(150:1);
		\draw (l) -- +(150:-1);
		\draw (l) -- +(30:-1);
		\node at ($0.5*(g)+0.5*(l)$) {$-$};
		
		\node (m) at (9, 0.5) {$\hat{v}$};
		\coordinate (n) at (9, -0.5);
		\draw [fill] (n) circle (0.05);
		\draw (m)--(n);
		\draw (m) -- +(30:1);
		\draw (m) -- ($(m)+(0,-1) +(30:-1)$);
		\draw (n) -- +(150:-1);
		\draw (n) -- ($(n)+(0,1)+(150:1)$);
		\node at ($0.5*(m)+0.5*(l)$) {$+$};
		
		\node (p) at (11, 0.5) {$\hat{v}$};
		\coordinate (q) at (11, -0.5);
		\draw [fill] (q) circle (0.05);
		\draw (p)--(q);
		\draw (p) -- +(30:1);
		\draw (p) -- +(150:1);
		\draw (q) -- +(150:-1);
		\draw (q) -- +(30:-1);
		\node at ($0.5*(m)+0.5*(q)$) {$=$};
		
		\node (r) at (13, -0.5) {$\hat{v}$};
		\coordinate (s) at (13, +0.5);
		\draw [fill] (s) circle (0.05);
		\draw (s)--(r);
		\draw (s) -- +(30:1);
		\draw (s) -- +(150:1);
		\draw (r) -- +(150:-1);
		\draw (r) -- +(30:-1);
		\node at ($0.5*(p)+0.5*(r)$) {$-$};
	\end{tikzpicture}
\end{figure*}

\noindent By repeatedly applying the same transformations to the right side, we obtain:

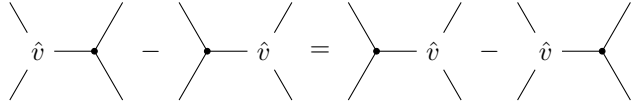
\begin{figure*}[!h]
	\centering
	\begin{tikzpicture}[scale=0.75]
		\node (a) at (-3, 0) {$\hat{v}$};
		\coordinate (b) at (-2, 0);
		\draw [fill] (b) circle (0.05);
		\draw (a) -- (b);		
		\draw (a) -- ($(a)+(120:1)$);
		\draw (a) -- ($(a)+(240:1)$);	
		\draw (b) -- ($(b)+(120:-1)$);
		\draw (b) -- ($(b)+(240:-1)$);	
		\node (d) at (1., 0) {$\hat{v}$};
		\coordinate (c) at (0., 0);
		\draw (c)--(d);
		\draw [fill] (c) circle (0.05);
		\draw (c) -- ($(c)+(120:1)$);
		\draw (c) -- ($(c)+(240:1)$);	
		\draw (d) -- ($(d)+(120:-1)$);
		\draw (d) -- ($(d)+(240:-1)$);	
		\node[] at ($0.5*(b) + 0.5*(c)$) {$-$};

		\coordinate (aa) at (3, 0);
		\node (bb) at (4, 0) {$\hat{v}$};
		\draw [fill] (aa) circle (0.05);
		\draw (aa) -- (bb);		
		\draw (aa) -- ($(aa)+(120:1)$);
		\draw (aa) -- ($(aa)+(240:1)$);	
		\draw (bb) -- ($(bb)+(120:-1)$);
		\draw (bb) -- ($(bb)+(240:-1)$);	
		\coordinate (dd) at (7, 0) ;
		\node (cc) at (6., 0) {$\hat{v}$};
		\draw (cc)--(dd);
		\draw [fill] (dd) circle (0.05);
		\draw (cc) -- ($(cc)+(120:1)$);
		\draw (cc) -- ($(cc)+(240:1)$);	
		\draw (dd) -- ($(dd)+(120:-1)$);
		\draw (dd) -- ($(dd)+(240:-1)$);	
		\node[] at ($0.5*(bb) + 0.5*(cc)$) {$-$};
		\node at ($0.5*(d)+0.5*(aa)$) {$=$};
	\end{tikzpicture}
\end{figure*}

This relation immediately implies the condition for the correctness of the multiplication in $\Lambda$ in Fig.\,\ref{fig:2diagrams}. It is important to note that this way one can multiply any connected diagram modulo the AS and IHX relations by $\Lambda$. Such diagrams include primitive Jacobi diagrams starting from the second order $\mathcal{P}_{n\geq 2}$ (at the zeroth and first orders, there are no internal vertices where the element of $\Lambda$ can be inserted). Hence, $\mathcal{P}_{\geq2}$ obtain the structure of $\Lambda$~-module \cite{vogel2011algebraic}.

The multiplication operation in the algebra $\Lambda$ is commutative. The proof of this fact can be found in the work of Vogel \cite[1301]{vogel2011algebraic}. Some elements of the $\Lambda$ algebra are shown in Fig.\,\ref{fig:some-of-Lambda}. A proof of the $\hat{x}_n$ diagrams to belong to the $\Lambda$ algebra can be found in \cite[1305]{vogel2011algebraic}.

\begin{figure}[!h]
	\centering
	\begin{tikzpicture}
		
		\node at (-5.7, 0.2) {$\hat{1}=$};
		\draw (-4.7,0)--+(90:1);
		\draw (-4.7,0)--+(-30:1);
		\draw (-4.7,0)--+(210:1);
		\draw [fill] (-4.7,0) circle (0.05);
		
		\node at (-1, 0.2) {$\hat{t}=$};
		\draw (90:0.3)--(90:1);
		\draw (-30:0.3)--(-30:1);
		\draw (210:0.3)--(210:1);
		\draw (0,0) circle (0.3);
		\draw [fill] (90:0.3) circle (0.05);
		\draw [fill] (-30:0.3) circle (0.05);
		\draw [fill] (210:0.3) circle (0.05);
		
		\node at (3.7, 0.2) {$\hat{x}_n=$};
		\draw (4., -0.5) -- (5.4, -0.5);
		\draw (4.7, 0.6) -- +(0, 0.4);
		\draw (4.7, 0.3) circle (0.3);
		\draw ($(4.7, 0.3)+(50:-0.3)$) -- (4.2, -0.5);
		\draw ($(4.7, 0.3)+(60:-0.3)$) -- (4.3, -0.5);
		\draw ($(4.7, 0.3)+(70:-0.3)$) -- (4.4, -0.5);
		\draw ($(4.7, 0.3)+(80:-0.3)$) -- (4.5, -0.5);
		\node at (4.7, -0.3) {$\cdot$};
		\node at (4.8, -0.3) {$\cdot$};
		\node at (4.9, -0.3) {$\cdot$};
		\draw ($(4.7, 0.3)+(120:-0.3)$) -- (5.1, -0.5);
		\draw ($(4.7, 0.3)+(130:-0.3)$) -- (5.2, -0.5);
	\end{tikzpicture} 
	\caption{Some elements of $\Lambda$.}
	\label{fig:some-of-Lambda}
\end{figure}

Despite many findings of Pierre Vogel there are still many features of $\Lambda$ that remain undiscovered. For example, dimensions of $\Lambda$ are only known up to degree $10$ \cite[\S 7.2]{chmutov2012introduction}.

Nevertheless, there exists a homomorphism $\varphi : \mathbb{Q}[t]\oplus\omega\mathbb{Q}[t,\sigma,\omega]\rightarrow\Lambda$ from the polynomial algebra $\mathbb{Q}[t]\oplus\omega\mathbb{Q}[t, \sigma,\omega]$ to the algebra $\Lambda$ (over the field of rational numbers $\mathbb{Q}$). The derivation of this homomorphism can be found in \cite[Section 5]{vogel2011algebraic}. The diagrams $\hat{x}_n\in\Lambda$ lie in the image of this homomorphism and can be expressed in terms of polynomials in $t$, $\sigma$ and $\omega$. In particular, 
\begin{equation}\label{x-n}
	\hat{x}_1 = 2 \hat{t} =2\varphi(t),\quad \hat{x}_2 = \hat{t}^2=\varphi(t^2),\quad \hat{x}_3 = \varphi (4 t^3 - \frac32 \omega),\quad \hat{x}_5 = \varphi (12 t^5 - \frac{17}2 t^2 \omega + \frac32 \sigma \omega).
\end{equation}

\begin{remark}
\label{remark32}
	The variable $t$ is of degree $1$ and corresponds to the diagram $\hat{t} = \varphi(t)$ shown in Fig.\,\ref{fig:some-of-Lambda}. The variable $\omega$~ is of degree $3$, and the combination of diagrams $\hat{\omega}_0 = \frac83 \hat{t}^3 - \frac23\hat{x}_3$ corresponds to it in $\Lambda$. The variable $\sigma$ of degree $2$ is missing in ``pure form'' and must always be in the product with $\omega$. In other words, there is no $\sigma$ diagram in $\Lambda$, but there are diagrams $\hat{\omega}_p=\varphi(\sigma^p\omega$) that are known to satisfy $\hat{\omega}_p\hat{\omega}_q =\hat{\omega}_0\hat{\omega}_{p+q}$, see \cite[Theorem 5.6.]{vogel2011algebraic}.
\end{remark}

Variables $t$, $\sigma$ and $\omega$ come with the package of marked diagrams described in Appendix~\ref{appendix}. 

It is known about this homomorphism $\varphi$ \cite{vogel1999universal} that at orders $\leq 10$ it is a bijection. It was hypothesized in \cite{vogel1999universal, kneissler2001spaces} that $\varphi$ is an isomorphism. However, later Pierre Vogel has showed \cite[Theorem 8.4.]{vogel2011algebraic} that in $\Lambda$ there is an element $\hat{P}$ of order $15$, such that $\hat{t}\hat{P} = 0$. Thus, there are zero divisors in $\Lambda$ and the homomorphism $\varphi$ has a nontrivial kernel.

\subsection{Characters on $\Lambda$}

It is the moment to return to the weight systems of Lie algebras discussed in Section \ref{sec:Lie-weights}. Their design can be extended to fixed diagrams with any number of legs, modulo AS and IHX relations. To do this, we put structure constants ${f^{abc}}$ at vertices and the edges are associated with $g_{ab}$.

For simplicity, we denote this generalized weight system for a Lie algebra $L$ as $\Phi_L$. $\Phi_L$ maps a fixed diagram with $n$ legs to a rank-$n$ tensor obtained from contracting the structure constants. To illustrate this, consider the diagram $\Psi$ shown in Fig.\,\ref{fig:psi-diag}. $\Phi_L(\Psi)=g_{ij}{f^{ia}}_c{f^{jb}}_d=(\Psi_L)^{ab}_{cd}$ is a rank-4 tensor that corresponds to the endomorphism of $L^{\otimes2}$.

\begin{figure}[h]
	\centering	
	\begin{tikzpicture}[scale=0.7]
		\node (a) at (-1.5, 1) {$a$};
		\node (b) at (-1.5, 0) {$b$};
		\node (c) at (1.5, 1) {$c$};
		\node (d) at (1.5, 0) {$d$};
		\coordinate (i) at (0, 1);
		\coordinate (j) at (0, 0);
		\draw [fill] (i) circle (0.05);
		\draw [fill] (j) circle (0.05);
		\draw (i)--(j)-- (d);
		\draw (i)-- (a);
		\draw (i)-- (c);
		\draw (j)-- (b);
	\end{tikzpicture}
	\caption{Diagram $\Psi$.}
	\label{fig:psi-diag}
\end{figure}

Now, in the case where $\hat{v}$ is an element of $\Lambda$, $\Phi_L(\hat{v})$ is a tensor of rank 3 obtained from the contraction of structure constants leaving three free indices. Since there are no other primitive tensors of rank 3 in a Lie algebra except for the structure constant $f$, it can be expected that $\Phi_L(\hat{v}) \propto f$. The proportionality coefficient is called the {\it character} $\chi_L(\hat{v})$ of an element $\hat{v}$. Then, for $\hat{v} \in \Lambda$ and $u$ being a diagram modulo AS and IHX relations, the equality $\Phi_L(\hat{v} u) = \chi_L(\hat{v})\Phi_L(u)$ holds. A rigorous and formal proof of this statement can be found in Vogel's paper \cite[Theorem 6.1]{vogel2011algebraic}.

In order to find characters of elements of the $\Lambda$-algebra of the form $\hat{t}^m \hat{\omega}_0^n \hat{\omega}_p$, we only need to know the values of $\chi_L(\hat t)$, $\chi_L(\hat\omega_0)$, and $\frac{\chi_L(\hat \omega_{p+1})}{\chi_L (\hat{ \omega}_p)}$, as this relation does not depend on $p$ since $\hat {\omega}_p\hat{ \omega}_q=\hat{\omega}_0 \hat{ \omega}_{p + q}$. These numbers are denoted by $t_L$, $\omega_L$ and $\sigma_L$ correspondingly. They depend on $L$ and a normalization of the metric tensor. We can say that they parameterize simple Lie algebras.

\begin{remark}
	In Vogel's paper \cite{vogel2011algebraic} there was no distinction between $t$, $ \hat{t} $ and $ t_L $. In this paper we use $\hat{v}$ when we refer to an element of $\Lambda$, i.e. a 3-legged diagram modulo the AS and IHX relation, and $v_L$ stands for a numerical value $\chi_L(\hat{v})$ with $L$ being a simple Lie algebra.  
\end{remark}

There is also a more convenient parameterization using three numbers $\alpha$, $\beta$ and $\gamma$. The construction of this parameterization is described in Sections 6 and 7 of Vogel's paper \cite{vogel2011algebraic}. Let us take a closer look at this construction. The key concept behind parameterization through $\alpha$, $\beta$, and $\gamma$ is a cubic relation for the operator $\Psi_L|_Y$ which is induced by $\Psi_L$. This operator acts on a space $Y$ defined as $S^2 L/\Omega$, where $\Omega$ is a one-dimensional space generated by the Casimir element $g_{ij}X^i\otimes X^j$ with $X^i$ being some basis in $L$.

The cubic relation is given by:

\begin{equation}
	(\Psi_L|_Y)^3 - t_L (\Psi_L|_Y)^2 + (\sigma_L-2t_L^2)\Psi_L|_Y- (\omega_L - t_L\sigma_L) = 0\,.
	\label{eq:cubic_equation}
\end{equation}
It was derived in Vogel's paper \cite{vogel2011algebraic} and can be used to decompose the space $Y$ into three subspaces corresponding to the eigenvalues $\alpha_L, \beta_L, \gamma_L$. These subspaces are denoted as $Y(\alpha)$, $Y(\beta)$, and $Y(\gamma)$. The three variables $\alpha, \beta$ and $\gamma$ are related to the variables $t, \sigma$ and $\omega$ as follows:

%	 Оn the one hand, the cubic equation can be written as equation (\ref{eq:cubic_equation}). On the other hand, it should look like $(\Psi_L|_Y - \alpha)(\Psi_L|_Y-\beta)(\Psi_L|_Y-\gamma) = 0$. By equating the corresponding coefficients, we get that 
\begin{equation}
	\alpha + \beta +\gamma = t, \quad \alpha \beta + \beta \gamma + \alpha \gamma = \sigma-2t^2, \quad \alpha \beta \gamma = \omega-t\sigma\,.
\end{equation}
The values $\alpha_L,\ \beta_L,\ \gamma_L$ for simple Lie algebras were calculated by Pierre Vogel in \cite{vogel2011algebraic} and are shown in Table \ref{tab:char} in the same normalization as in \cite{mkrtchyan2012casimir}.

\begin{table}[h]
	\centering
	\begin{tabular}{|c|c|c|c|}
		\hline
		Lie algebra $L$& $\alpha_L$ & $\beta_L$ & $\gamma_L$ \\
		\hline
		$\mathfrak{sl}_{n}$ & $-2$ & $2$ &$n$ \\
		\hline
		$\mathfrak{so}_{n}$ & $-2$ & $4$ & $n-4$\\
		\hline
		$\mathfrak{sp}_{2n}$ & $-2$ & $1$ &$n+2$ \\
		\hline
		$G_2$ & $-2$ & $10/3$ & $8/3$ \\
		\hline
		$F_4$ & $-2$ & $5$ &$6$ \\
		\hline
		$E_6$ & $-2$ & $6$ & $8$ \\
		\hline
		$E_7$ & $-2$ & $8$ & $12$\\
		\hline
		$E_8$ & $-2$ & $12$ & $20$ \\
		\hline
	\end{tabular}
	\caption{Vogel parameters.}	
	\label{tab:char}
\end{table}

Not only the characters on $\Lambda$, but also some other properties of Lie algebras can be expressed as functions of $\alpha$, $\beta$ and $\gamma$.  This concept is referred to as \textit{Vogel universality}. For example, it is possible to express the dimension of a Lie algebra by the following formula \cite{vogel2011algebraic}:

\begin{equation}
	\dim L = \frac{(\alpha_L - 2t_L)(\beta_L - 2t_L)(\gamma_L - 2t_L)}{\alpha_L \beta_L \gamma_L }\,.
	\label{eq:universal-dim}
\end{equation}

There are other features that can be universalized:

\begin{enumerate}
	\item[1)] dimensions of representations, appearing in decomposition of power of adjoint representation~\cite{landsberg2006universal}; 
	\item[2)] eigenvalues of higher Casimir operators~\cite{mkrtchyan2012casimir,isaev2022split}; 
	\item[3)] volume of simple Lie groups~\cite{khudaverdian2017universal}; 
	\item[4)] quantum knot invariants~\cite{mironov2016universal}; etc.
\end{enumerate}

Notably, all known universal formulae work only in the adjoint sector, i.e. for tensor powers of the adjoint representation. It was hypothesized that there exists some object that generalizes Lie algebras. It was referred to as \textit{universal Lie algebra}. Some attempts to describe it were made by Vogel in his unpublished works \cite{vogel1999universal}. Currently, it remains unknown whether it exists or not.

\setcounter{equation}{0}
\section{Construction of weight system kernel}\label{sec:weight-system-kernel}

In this section, we describe a method for construction of the kernel of Lie algebra weight system and apply it to provide lowest order diagrams from $\varphi_{\mathfrak{sl}_n}$ kernel.

\subsection{Method}
\label{method}

Our method to construct the kernel of a Lie algebra weight system is based on the technique introduced by Pierre Vogel to construct a 17-order Jacobi diagram in the kernel of all Lie algebra weight systems (see \cite{vogel2011algebraic}). The underlying idea is to multiply some Jacobi diagram by a specific element of $\Lambda$ that is killed by the corresponding character.

To construct such diagram, one should remember from the previous section that elements of $\Lambda$ can be expressed as symmetric polynomials in $\alpha$, $\beta$, $\gamma$. The Lie algebra character $\chi_L$ maps this polynomial to its value on Vogel parameters $\alpha_L$, $\beta_L$, $\gamma_L$ that can be found in Table \ref{tab:char}.

It can be seen from Table \ref{tab:char} that for Lie algebra series $\mathfrak{sl}_n$ there is a symmetric polynomial $P_{sl} = (\alpha+\beta)(\beta + \gamma)(\alpha+\gamma) = 2 t \sigma - \omega - 2 t^3$ that is killed by character $\chi_{sl}$. Similar polynomials exist for $\mathfrak{so}_n$ series and for the exceptional Lie algebras \cite{vogel2011algebraic}. These are polynomials in $\mathbb{Q}[t, \sigma, \omega]$ that are not in $\mathbb{Q}[t] \oplus \omega \mathbb{Q}[t, \sigma, \omega]$. Hence, they do not correspond to any element of $\Lambda$ as they all contain $\sigma$ in pure form. 

Our method is a modification that works for a polynomial $P$ in $\mathbb{Q}[t, \sigma, \omega]$ that is not in $\mathbb{Q}[t]\oplus \omega \mathbb{Q}[t, \sigma, \omega]$. We provide several ways to multiply $P$ on Jacobi diagrams to create diagrams in the kernel of a specific Lie algebra weight system.

\begin{enumerate}
	\item Multiply any primitive Jacobi diagram of order $\geq 2$ (in the $\Lambda$-algebra sense) by an element $\widehat{(\omega P)} \in \Lambda$ or by the multiples of it.
	
	\item Multiply marked diagrams described in Appendix \ref{appendix} by $P$ with multiplication by $\sigma$ acting as raising of the total marking.
	
	\item Multiply special linear combinations of diagrams given in Appendix \ref{appendix} in Fig.\,\ref{fig:two_additional_relations} by $P$. These combinations of diagrams allow for multiplication by $\sigma$ only when considered under an action of Lie algebra weight system.
	\end{enumerate}
	
	%To find all diagrams in the kernel, o
	One also needs to check whether the resulting diagrams are linearly independent. To check linear independence, it is very useful to operate in open Jacobi diagrams, because they have an extra grading by the number of legs. Another way to check the linear independence of diagrams is to calculate some weight system that can distinguish these diagrams.

\subsection{Example for $\mathfrak{sl}_n$ weight system}

The estimates for dimensions of the image of $\mathfrak{sl}_n$ weight system have been found in \cite{lanina2021chern,lanina2022implications}. Dimensions of $\mathcal{P}_n$ and $\mathcal{A}_n$ were calculated in \cite{vassiliev1990homological,bar1995vassiliev,kneissler1997number}. Subtraction yields the estimates for the dimensions of $\text{Ker}~\varphi_{\mathfrak{sl}_n}$. All these estimates are listed in Table \ref{tab:dimensions}. 

	\begin{table}[!h]
	\centering	
	\begin{tabular}{|c|c|c|c|c|c|c|c|c|c|c|c|c|c|}
		\hline
		\raisebox{-0.1cm}{$n$} & \raisebox{-0.1cm}{1} & \raisebox{-0.1cm}{2} & \raisebox{-0.1cm}{3} & \raisebox{-0.1cm}{4} & \raisebox{-0.1cm}{5} & \raisebox{-0.1cm}{6} & \raisebox{-0.1cm}{7} & \raisebox{-0.1cm}{8} & \raisebox{-0.1cm}{9} & \raisebox{-0.1cm}{10} & \raisebox{-0.1cm}{11} & \raisebox{-0.1cm}{12} & \raisebox{-0.1cm}{13} \\ [1.5ex]
		\hline
		\hline
		\raisebox{-0.1cm}{$\dim \mathcal{P}_n\; =$} & \raisebox{-0.1cm}{0} & \raisebox{-0.1cm}{1} & \raisebox{-0.1cm}{1} & \raisebox{-0.1cm}{2} & \raisebox{-0.1cm}{3} & \raisebox{-0.1cm}{5} & \raisebox{-0.1cm}{8} & \raisebox{-0.1cm}{12} & \raisebox{-0.1cm}{18} & \raisebox{-0.1cm}{27} & \raisebox{-0.1cm}{39} & \raisebox{-0.1cm}{55} & \raisebox{-0.1cm}{?} \\ [1.5ex]
		\hline
		\raisebox{-0.1cm}{$\dim \varphi_{\mathfrak{sl}_n} (\mathcal{P}_n) \; \geq$} & \raisebox{-0.1cm}{0} & \raisebox{-0.1cm}{1} & \raisebox{-0.1cm}{1} & \raisebox{-0.1cm}{2} & \raisebox{-0.1cm}{3} & \raisebox{-0.1cm}{5} & \raisebox{-0.1cm}{8} & \raisebox{-0.1cm}{11} & \raisebox{-0.1cm}{16} & \raisebox{-0.1cm}{22} & \raisebox{-0.1cm}{30} & \raisebox{-0.1cm}{42} & \raisebox{-0.1cm}{53} \\ [1.5ex]
		\hline
		\raisebox{-0.1cm}{$\dim \text{Ker}\, \varphi_{\mathfrak{sl}_n} (\mathcal{P}_n) \; \leq$} & \raisebox{-0.1cm}{\textbf{0}} & \raisebox{-0.1cm}{\textbf{0}} & \raisebox{-0.1cm}{\textbf{0}} & \raisebox{-0.1cm}{\textbf{0}} & \raisebox{-0.1cm}{\textbf{0}} & \raisebox{-0.1cm}{\textbf{0}} & \raisebox{-0.1cm}{\textbf{0}} & \raisebox{-0.1cm}{\textbf{1}} & \raisebox{-0.1cm}{\textbf{2}} & \raisebox{-0.1cm}{\textbf{5}} & \raisebox{-0.1cm}{\textbf{9}} & \raisebox{-0.1cm}{\textbf{13}} & \raisebox{-0.1cm}{?} \\ [1.5ex]
		\hline \hline
		\raisebox{-0.1cm}{$\dim \mathcal{A}_n\; = $} & \raisebox{-0.1cm}{0} & \raisebox{-0.1cm}{1} &\raisebox{-0.1cm}{1} & \raisebox{-0.1cm}{3} & \raisebox{-0.1cm}{4} & \raisebox{-0.1cm}{9} & \raisebox{-0.1cm}{14} & \raisebox{-0.1cm}{27} & \raisebox{-0.1cm}{44} & \raisebox{-0.1cm}{80} & \raisebox{-0.1cm}{132} & \raisebox{-0.1cm}{232} & \raisebox{-0.1cm}{?} \\ [1.5ex]
		\hline
		\raisebox{-0.1cm}{$\dim \varphi_{\mathfrak{sl}_n} (\mathcal{A}_n) \; \geq$} & \raisebox{-0.1cm}{0} & \raisebox{-0.1cm}{1} & \raisebox{-0.1cm}{1} & \raisebox{-0.1cm}{3} & \raisebox{-0.1cm}{4} & \raisebox{-0.1cm}{8} & \raisebox{-0.1cm}{11} & \raisebox{-0.1cm}{19} & \raisebox{-0.1cm}{25} & \raisebox{-0.1cm}{39} & \raisebox{-0.1cm}{50} & \raisebox{-0.1cm}{75} & \raisebox{-0.1cm}{95} \\ [1.5ex]
		\hline 
		\raisebox{-0.1cm}{$\dim \text{Ker}\, \varphi_{\mathfrak{sl}_n} (\mathcal{A}_n) \; \leq $} & \raisebox{-0.1cm}{\textbf{0}} & \raisebox{-0.1cm}{\textbf{0}} & \raisebox{-0.1cm}{\textbf{0}} & \raisebox{-0.1cm}{\textbf{0}} & \raisebox{-0.1cm}{\textbf{0}} & \raisebox{-0.1cm}{\textbf{1}} & \raisebox{-0.1cm}{\textbf{3}} & \raisebox{-0.1cm}{\textbf{8}} & \raisebox{-0.1cm}{\textbf{19}} & \raisebox{-0.1cm}{\textbf{41}} & \raisebox{-0.1cm}{\textbf{82}} & \raisebox{-0.1cm}{\textbf{157}} & \raisebox{-0.1cm}{?} \\ [1.5ex]
		\hline
	\end{tabular}	
	\caption{Estimates for the dimensions of the kernel of $\varphi_{\mathfrak{sl}_n}.$}
	\label{tab:dimensions}
\end{table}

In this subsection, the kernel of $\mathfrak{\varphi_{\mathfrak{sl}_n}}$ at first few orders in both primitive and non-primitive cases are constructed using the described method. It is demonstrated below that the inequalities in Table \ref{tab:dimensions} at the lowest orders turn out to be nicely saturated, what also validates the results of our previous study of \cite{lanina2021chern,lanina2022implications}.

\subsubsection{Primitive diagrams in kernel}\label{sec:prim-ker-sl}

For the $\mathfrak{sl}_N$ weight system we have polynomial $P_{sl} = (\alpha+\beta)(\beta + \gamma)(\alpha+\gamma) = 2 t \sigma - \omega - 2 t^3$ that vanishes when taken on values of the parameters corresponding to $\mathfrak{sl}_n$ Lie algebra, see Table \ref{tab:char}.
%is killed by the extension of character $\chi_{sl}$. 

$\widehat{(\omega P_{sl})}$ is an element of $\Lambda$ that can be expressed in terms of $\hat{x}_n$ using~\eqref{x-n} as $\frac43 \hat{t} \hat{x}_5 + \frac{16}9 \hat{t}^6-\frac83 \hat{t}^3 \hat{x}_3 - \frac49 \hat{x}_3^2 $, see Fig.\,\ref{fig:omPsl}. That element is killed by $\mathfrak{sl}_n$ weight system and it can be multiplied by any Jacobi diagram to create a diagram in the kernel of the $\mathfrak{sl}_n$ weight system.

	\begin{figure}[h]
	\centering
	\begin{tikzpicture}[scale=0.75]
		\node[scale=1.4] at (0.4, 0) {$\frac23$};

		\xVogel{(1, 1)}{(1, -1)}{(2.6, 0)}{5};
		\bubble{(2.5,0)}{(3.41, 0)};
		\node[scale=1.4] at (4.2, 0) {$+ \frac{1}{36}$};

		\coordinate (x) at (5,1);
		\coordinate (y) at (5, -1);
		\coordinate (c) at (6,0);
		\coordinate (z) at (7.41, 0);

		\bubble{($1.*(x)+0.*(c) $)}{($0.5*(x)+0.5*(c) $)};
		\bubble{($0.5*(x)+0.5*(c) $)}{($0.*(x)+1*(c) $)};
		
		\bubble{($1.*(y)+0.*(c) $)}{($0.5*(y)+0.5*(c) $)};
		\bubble{($0.5*(y)+0.5*(c) $)}{($0.*(y)+1*(c) $)};
		
		\bubble{($1.*(z)+0.*(c) $)}{($0.5*(z)+0.5*(c) $)};
		\bubble{($0.5*(z)+0.5*(c) $)}{($0.*(z)+1*(c) $)};

		\node[scale=1.4] at (8.2, 0) {$- \frac{1}{3}$};

		\xVogel{(9, 1)}{(9,-1)}{(10.4,0)}{3};
		\bubble{(10.4, 0)}{(10.8, 0)};
		\bubble{(10.7, 0)}{(11.1, 0)};
		\bubble{(11.0, 0)}{(11.41, 0)};
		\node[scale=1.4] at (12.2, 0) {$- \frac{4}{9}$};

		\xVogel{(13, 1)}{(13, -1)}{(15.41,0)}{3};
		\xVogel{(14, 0.4)}{(14, -0.4)}{(15,0)}{3};
	\end{tikzpicture}	
	\caption{$\widehat{\omega P_{sl}}$~--- the first element of $\Lambda$ that is killed by $\chi_{sl}$.}
	\label{fig:omPsl}
\end{figure}

The first such diagram can be obtained by taking the product $\widehat{(\omega P_{sl})}\bbl$.  Other diagrams are depicted in Table \ref{tab:primitive}. Note that they are all open Jacobi diagrams that have an extra grading by the number of legs $l$. Relation between open and closed Jacobi diagrams is described in Section \ref{sec:open_Jacobi}.

	\begin{table}[!h]
	\centering
	\begin{tabular}{ c| >{\centering\arraybackslash}m{3cm} |>{\centering\arraybackslash}m{3cm}|>{\centering\arraybackslash}m{3cm}| }
		\backslashbox{$l$}{$n$} & $8$ & $9$ & $10$
		\tabularnewline \hline 
		2
		&
		\scalebox{0.8}{\begin{tikzpicture}
				\bubble{(0,0)}{(2,0)};
				\node at (-1, 0) { $ (\widehat{\omega P_{sl}})$};
		\end{tikzpicture}} & \scalebox{0.8}{\begin{tikzpicture}
				\bubble{(0,0)}{(2,0)};
				\node at (-1, 0) { $\hat{t} (\widehat{\omega P_{sl}})$};
		\end{tikzpicture}} 	& \vspace{1mm}
		\raisebox{0.15cm}{\scalebox{0.8}{\begin{tikzpicture}
				\bubble{(0,0)}{(2,0)};
				\node at (-1, 0) { $\hat{t}^2 (\widehat{ \omega P_{sl}})$};
		\end{tikzpicture}}}  \scalebox{0.8}{\begin{tikzpicture}
				\bubble{(0,0)}{(2,0)};
				\node at (-1, 0) { $ (\widehat{\sigma \omega P_{sl}})$};
		\end{tikzpicture}}
		\tabularnewline \hline 
		4 
		&
		&
		
		\scalebox{0.7}{
			\begin{tikzpicture}
				\coordinate (a) at (0,0);
				\draw (a) -- + (1,-0.5);
				\draw (a) -- + (1,0.5);
				\draw (a) -- + (-1,0.5);
				\draw (a) -- + (-1,-0.5);
				\node[scale=0.7] at ($(a)+(0, -0.3)$) {$0$};
				\node[scale=1.5] at ($(a)+(-1.5,0)$) {$P_{sl}$};
				\draw[fill] (a) circle (0.1);
			\end{tikzpicture}
		}
		& \vspace{1mm}
		\raisebox{0.2cm}{\scalebox{0.7}{
			\begin{tikzpicture}
				\coordinate (a) at (0,0);
				\draw (a) -- + (1,-0.5);
				\draw (a) -- + (1,0.5);
				\draw (a) -- + (-1,0.5);
				\draw (a) -- + (-1,-0.5);
				\node[scale=0.7] at ($(a)+(0, -0.3)$) {$0$};
				\node[scale=1.5] at ($(a)+(-1.5,0)$) {$\hat{t}P_{sl}$};
				\draw[fill] (a) circle (0.1);
			\end{tikzpicture}
		}}
		\newline 
		\scalebox{0.7}{\begin{tikzpicture}
				\coordinate (d) at (0,0);
				\draw (d) +(-1, 0.5) -- + (1, 0.5);
				\draw (d) +(-1, -0.5) -- + (1, -0.5);
				\draw (d) +(0.3, -0.5) -- + (0.3, 0.5);
				\draw (d) +(-0.3, -0.5) -- + (-0.3, 0.5);
				\node at ($(d)+(-1.6,0)$) {($\widehat{\omega P_{sl}}$)};
			\end{tikzpicture}
		}
		\tabularnewline
		\hline
		6 & & & \vspace{1mm}
		\raisebox{0.1cm}{\scalebox{0.4}{
		\begin{tikzpicture}
			\coordinate (a) at (0,0);
			\coordinate (b) at (3, 0);
			\node[scale=2] at ($0.5*(a)+0.5*(b)$) {$-\hat{t}$};
			\node[scale=2] at ($(a)+(-1.7, 0)$) {$P_{sl}($};
			\node[scale=2] at ($(b)+(1.2, 0)$) {$)$};
		\draw (a) circle (0.5);
		\draw (a) +(30:1) -- +(30:0.5);
		\draw (a) +(90:1) -- +(90:0.5);
		\draw (a) +(330:1) -- +(330:0.5);
		\draw (a) +(150:1) -- +(150:0.5);
		\draw (a) +(210:1) -- +(210:0.5);
		\draw (a) +(270:1) -- +(270:0.5);
			\draw (a) + (-0.5, 0)-- +(0.5, 0);

			\draw (b) circle (0.5);
			\draw (b) +(30:1) -- +(30:0.5);
			\draw (b) +(90:1) -- +(90:0.5);
			\draw (b) +(330:1) -- +(330:0.5);
			\draw (b) +(150:1) -- +(150:0.5);
			\draw (b) +(210:1) -- +(210:0.5);
			\draw (b) +(270:1) -- +(270:0.5);
			
		\end{tikzpicture}
		}}
	\end{tabular}
	\caption{Diagrams in the kernel of $\varphi_{\mathfrak{sl}_n}$.}
	\label{tab:primitive}
\end{table}

The 6-legged diagram in Table \ref{tab:primitive} can be multiplied by $\sigma$ because of the top formula in Fig.\,\ref{fig:two_additional_relations}. The action of $\sigma$ on this diagram is shown in Fig.\,\ref{fig:6legs10ord}.

\begin{figure}[h!]
	\centering	
	\begin{tikzpicture}[scale=0.7]
			\coordinate (a) at (-1.0,0);
		\coordinate (b) at (2.5, 0);
		\node[scale=1.2] at ($0.8*(a)+0.5*(b)$) {$-$};
		\node[scale=1.5] at ($0.5*(a)+0.65*(b)$) {$\hat{t}$};
		\node[scale=1.5] at ($(a)+(-1.8, 0)$) {$\sigma$};
		\node[scale=2.2] at ($(a)+(-1.2, 0)$) {$($};
		\node[scale=2.2] at ($(b)+(1.2, 0)$) {$)$};
		\draw (a) circle (0.5);
		\draw (a) +(30:1) -- +(30:0.5);
		\draw (a) +(90:1) -- +(90:0.5);
		\draw (a) +(330:1) -- +(330:0.5);
		\draw (a) +(150:1) -- +(150:0.5);
		\draw (a) +(210:1) -- +(210:0.5);
		\draw (a) +(270:1) -- +(270:0.5);
		\draw (a) + (-0.5, 0)-- +(0.5, 0);

		\draw (b) circle (0.5);
		\draw (b) +(30:1) -- +(30:0.5);
		\draw (b) +(90:1) -- +(90:0.5);
		\draw (b) +(330:1) -- +(330:0.5);
		\draw (b) +(150:1) -- +(150:0.5);
		\draw (b) +(210:1) -- +(210:0.5);
		\draw (b) +(270:1) -- +(270:0.5);
		\node[scale=1.2] at (4.5, 0) {$=$};
			\coordinate (c) at (6,0);
		\coordinate (d) at (10, 0);
		\node[scale=1.2] at ($0.5*(c)+0.5*(d)+(-0.45, 0)$) {$-$};
		\node[scale=1.5] at ($0.5*(c)+0.5*(d)+(0.45, 0)$) {$\hat{\omega}$};
		\draw (c) circle (0.5);
		\draw (c) +(30:1) -- +(30:0.5);
		\draw (c) +(90:1) -- +(90:0.5);
		\draw (c) +(330:1) -- +(330:0.5);
		\draw (c) +(150:1) -- +(150:0.5);
		\draw (c) +(210:1) -- +(210:0.5);
		\draw (c) +(270:1) -- +(270:0.5);
		\draw (c) + (-0.5, 0)-- +(0.5, 0);
		\draw[fill] (c) circle(0.07);
		\node[scale=0.75] at ($(c)+(0, -0.25)$) {$0$};

		\draw (d) circle (0.5);
		\draw (d) +(30:1) -- +(30:0.5);
		\draw (d) +(90:1) -- +(90:0.5);
		\draw (d) +(330:1) -- +(330:0.5);
		\draw (d) +(150:1) -- +(150:0.5);
		\draw (d) +(210:1) -- +(210:0.5);
		\draw (d) +(270:1) -- +(270:0.5);
		
	\end{tikzpicture}
	\caption{Action of $\sigma$ on the 6-legged diagram.}
	\label{fig:6legs10ord}
\end{figure}

Linear independence of the diagrams depicted in Table \ref{tab:primitive} can be checked by calculating the value of $\mathfrak{so}_n$ weight system on these diagrams. One can compare the number of these diagrams with the estimates in Table \ref{tab:dimensions} and see that the constructed diagrams fully generate $\text{Ker}~\varphi_{\mathfrak{sl}_n}$ among the primitive Jacobi diagrams up to the 10-th order.

\subsubsection{Non-primitive diagrams in kernel}

As can be seen from Table \ref{tab:dimensions}, non-primitive diagrams first appear in the kernel of $\mathfrak{sl}_n$ weight systems at the 6-th order. In paper \cite{lanina2022implications}, the element of $\mathcal{C}_6$ in the kernel of $\varphi_{\mathfrak{sl}_n}$ is written explicitly:

\begin{figure}[!h]
	\centering
	\begin{tikzpicture}[scale=0.64]
		\coordinate (a) at (-2, 0);
		\coordinate (b) at (1, 0);
		\node at ($0.5*(a) + 0.5*(b)$) {$-$};
		
		\draw (a) +(30:1)-- +(150:1);
		\draw (a) +(-30:1)-- +(-150:1);
		\draw (a) +(0.4, 0.5) -- +(60:1);
		\draw (a) +(-0.4, 0.5) -- +(120:1);
		\draw (a) + (0, 0.5)-- + (0, 1);
		\draw (a) + (0, -0.5)-- +(0, -1);

		\draw (b) +(30:1)-- +(150:1);
		\draw (b) +(-30:1)-- +(-150:1);
		\draw (b) +(0.3, 0.5) -- +(70:1);
		\draw (b) +(-0.3, 0.5) -- +(110:1);
		\draw (b) + (0.3, -0.5)-- + (110:-1);
		\draw (b) + (-0.3, -0.5)-- +(70:-1);

		\draw[thick] (a) circle (1);
		\draw[thick] (b) circle (1);
	\end{tikzpicture}
	\caption{Non-primitive element of $\text{Ker}~\varphi_{\mathfrak{sl}_n}$ at the 6-th order.}
	\label{fig:6th_order_nonprim}
\end{figure} 

Making use of $\Lambda$, this combination of diagrams can be rewritten as $(\bbl)(\hat{t}^2 \bbl) - (\hat{t} \bbl)(\hat{t} \bbl)$. Now it is easy to see that inclusion of this diagram in the kernel of $\varphi_{\mathfrak{sl}_n}$ is explained by multiplicativity of Lie algebra weight systems. Indeed, for any Lie algebra weight system, $ \varphi_L\left( (\bbl)(\hat{t}^2 \bbl) \right) = \varphi_L(\bbl) \varphi_L(\hat{t}^2 \bbl) = t_L^2 \varphi_L(\bbl)^2 $. On the other side, $ \varphi_L \left( (\hat{t} \bbl)(\hat{t} \bbl)  \right) =t_L^2 \varphi_L(\bbl)^2$, and we get the desired equality.

The same way, one can construct other non-primitive elements of space of closed Jacobi diagrams $\mathcal{C}$ in $\text{Ker}~\varphi_{\mathfrak{sl}_n}$ at higher orders. At the seventh order there are two such diagrams, both are depicted in Fig.\,\ref{fig:7ord_nonprim}.

\begin{figure}[!h]
	\centering
	\begin{tikzpicture}[scale=0.64]
		\coordinate (a) at (-2, 0);
		\coordinate (b) at (1, 0);
		\node at ($0.5*(a) + 0.5*(b)$) {$-$};
		
		\draw (a) +(30:1)-- +(150:1);
		\draw (a) +(-30:1)-- +(-150:1);
		\draw (a) +(0.5, 0.5) -- +(50:1);
		\draw (a) +(-0.5, 0.5) -- +(130:1);
		\draw (a) +(0.16, 0.5) -- +(75:1);
		\draw (a) +(-0.16, 0.5) -- +(105:1);
		
		\draw (a) + (0, -0.5)-- +(0, -1);

		\draw (b) +(30:1)-- +(150:1);
		\draw (b) +(-30:1)-- +(-150:1);
		
		\draw (b) + (0.3, -0.5)-- + (110:-1);
		\draw (b) + (-0.3, -0.5)-- +(70:-1);
		
		\draw (b) +(0.4, 0.5) -- +(60:1);
		\draw (b) +(-0.4, 0.5) -- +(120:1);
		\draw (b) + (0, 0.5)-- + (0, 1);

		\draw[thick] (a) circle (1);
		\draw[thick] (b) circle (1);
		
		\coordinate (c) at (8, 0);
		
		\draw[thick] (c) circle (1);
		\draw (c)+(0,0.3) circle (0.3);
		\draw (c)+($(0,0.3)+(45:0.3)$)-- +(60:1);
		\draw (c)+($(0,0.3)+(135:0.3)$)--+(120:1);
		\draw (c)+($(0,0.3)+(225:0.3)$)--+(180:1);
		\draw (c)+($(0,0.3)+(315:0.3)$)--+(0:1);
		
		\draw ($(c)+(-30:1)$)--($(c)+(30:-1)$);
		\draw (c) +(0.3, -0.5)-- +(110:-1);
		\draw (c) +(-0.3, -0.5)-- +(70:-1);
		
		\coordinate (d) at (11, 0);
		
		\draw[thick] (d) circle (1);
		\draw (d)+(0,0.3) circle (0.3);
		\draw (d)+($(0,0.3)+(45:0.3)$)-- +(60:1);
		\draw (d)+($(0,0.3)+(135:0.3)$)--+(120:1);
		\draw (d)+($(0,0.3)+(225:0.3)$)--+(180:1);
		\draw (d)+($(0,0.3)+(315:0.3)$)--+(0:1);
		
		\draw ($(d)+(-30:1)$)--($(d)+(30:-1)$);
		\draw (d) +(0, -0.5)-- +(0,-1);
		\draw (d) +(30:1) -- + ($0.5*(0:1)+0.5*(0,0.3)+0.5*(315:0.3) $);
		\node at ($0.5*(c) + 0.5*(d)$) {$-$};
	\end{tikzpicture}
	\caption{Two elements of $\mathcal{C}_7$ in the kernel of $\varphi_{\mathfrak{sl}_n}. $}
	\label{fig:7ord_nonprim}
\end{figure}

Nevertheless, according to Table \ref{tab:dimensions}, at 7-th order there are at most 3 such linearly independent elements. The third diagram can be constructed using the bottom formula from Fig.\,\ref{fig:two_additional_relations}:

	\begin{figure}[!h]
	\centering
	\begin{tikzpicture}[scale=0.64]
		\coordinate (a) at (-7, 0);
		\coordinate (b) at (-3.4, 0);
		\coordinate (c) at (0, 0);
		\coordinate (d) at (4.7, 0);
		\node[scale=0.84] at ($0.5*(a)+0.5*(b)$) {$+ \frac43 \hat{\omega} \hat{t}$};
		\node[scale=0.84] at ($0.5*(c)+0.5*(b)$) {$+ 2 \hat{t}$};
		\node[scale=0.84] at ($0.5*(c)+0.5*(d)$) {$+ (\hat{t}^3 - \frac12 \hat{\omega})$};
		\node[scale=0.84] at ($(a) + (-2.3, 0)$) {$-(\hat{\omega}+2 \hat{t}^3)$};
		\draw[thick] (a) circle (1);
		\draw (a) circle (0.4);
		\draw (a) + (45:0.4)-- +(45:1);
		\draw (a) + (135:0.4)-- +(135:1);
		\draw (a) + (225:0.4)-- +(225:1);
		\draw (a) + (315:0.4)-- +(315:1);
		
		\draw[thick] (b) circle(1);
		\draw (b) +(30:1)-- + (150:1);
		\draw (b) +(30:-1)-- + (150:-1);
		\draw (b) +(0, 0.5)-- +(0, -0.5);
		
		\draw[thick] (c) circle(1);
		\draw (c) +(45:1)-- + (45:-1);
		\draw (c) +(135:1)-- + (135:-1);
		\draw[fill] (c) circle (0.1);
		\node[scale=0.8] at ($(c)+(0, -0.4)$) {$0$};
		
		\draw[thick] (d) circle (1);
		
		\bubble{($(d)+(25:1)$)}{($(d)+(155:1)$)};
		\bubble{($(d)+(25:-1)$)}{($(d)+(155:-1)$)};

	\end{tikzpicture}
	\caption{The third nonprimitive element of $\text{Ker}~\varphi_{\mathfrak{sl}_n}$ at 7-th order.}

\end{figure}

\setcounter{equation}{0}
\section{Relation between Chern-Simons invariants and Kontsevich integral}\label{sec:CS-KI}

In this section, we discuss the connection of the obtained results with the 3D quantum Chern-Simons field theory. Detailed description of the connection between the Chern-Simons theory and Vassiliev invariants can be found in particular in \cite{dunin2013kontsevich}.

The 3-dimensional Chern-Simons action with a simple compact Lie gauge group $G$ for the vector field $A_{\mu} = A_{\mu}^aT_a$, where $T_a$ are generators of the corresponding Lie algebra $L$, is given by:
\be
S[A] = \dfrac{\kappa}{4\pi}\int_{\mathbb{R}^3}d^3x \,\epsilon^{\mu\nu\rho}\,\Tr\left(A_{\mu}\partial_{\nu}A_{\rho} + \frac{2}{3}A_{\mu}A_{\nu}A_{\rho}\right).
\ee
Correlators in the Chern--Simons theory can be obtained from the gauge invariant Wilson loops (called the colored HOMFLY polynomials) by expansion of the path ordered exponent:
\begin{equation}
	\label{WilsonLoopExpValue}
\begin{aligned}
		\langle W_R(K)\rangle &=\frac{1}{\dim(R)}\left\langle \text{tr}_{R} \ \text{Pexp} \left( \oint_{K} A_\mu dx^\mu \right) \right\rangle_{\text{CS}}=1\, + \,\frac{\text{tr}_R\,(T_a)}{\dim(R)} \oint_{K} \left\langle A_\mu^a(x) \right\rangle dx^\mu\,+ \\
		&+\,\frac{\text{tr}_R\, (T_{a_1}T_{a_2})}{\text{qdim}(R)}\oint_{K} dx_2^{\mu_2} \int_0^{x_2}dx_1^{\mu_1} \left\langle A_{\mu_1}^{a_1}(x_1) A_{\mu_2}^{a_2}(x_2)\right\rangle\, +\, \dots=\sum_{n=0}^{\infty}\hbar^n \sum_{m=1}^{\dim\,\mathcal{A}_{n}} v_{n,m}(K)\,\mathcal{G}_{n, m}^R
\end{aligned}
\end{equation}
where 
\begin{equation}\label{V-inv}
	\begin{aligned}
		v_{n, m}(K) \sim \oint d x_{2n}^{\mu_{2n}} \ldots \int d x_{1}^{\mu_{1}}\left\langle A_{\mu_1}^{a^{(m)}_{1}}\left(x_{1}\right) \ldots  A_{\mu_{2n}}^{a^{(m)}_{2n}}\left(x_{2n}\right)\right\rangle , 
	\end{aligned}
\end{equation}
and
\begin{equation}\label{G-factors}
	\mathcal{G}_{n, m}^R \sim \text{tr}_{R} \left(T_{a^{(m)}_1} T_{a^{(m)}_2}\ldots T_{a^{(m)}_{2n}} \right).
\end{equation}

The answer~\eqref{WilsonLoopExpValue} is very similar to the Kontsevich integral~\eqref{ki}. The only difference is that in the Kontsevich integral instead of group factors $\mathcal{G}_{n, m}^R$, there are basic elements $D_{n,m}$ of the algebra of chord diagrams. Actually, this is not just a similarity. In fact, Lie algebras weight systems map chord diagrams to group factors of the Wilson loop expansion~\eqref{WilsonLoopExpValue}, and due to linearity of the mapping

\begin{equation}
	\langle W_R(K)\rangle = \varphi_L^R\left(I(K)\right).
\end{equation} 
Thus, the results on the kernel of Lie algebras weight systems acquire the following meaning for the Chern--Simons theory. In a fixed Lie algebra $L$, there appear relations between group factors $\mathcal{G}_{n, m}^R$ at a fixed order $n$. Therefore, not all correlators $v_{n,m}(K)$ can be extracted from the Wilson loop~\eqref{WilsonLoopExpValue} for the Chern--Simons theory with a fixed gauge group. In particular, one primitive correlator is absent already at the 8-th order in the case of $SU(N)$ gauge group\footnote{Wilson loops in the 3D Chern--Simons theory with $SU(N)$ gauge group are quantum knot invariants called colored (by a representation $R$ of $\mathfrak{sl}_n$ Lie algebra) HOMFLY polynomials.}, see Section~\ref{sec:prim-ker-sl}.

 \section{On knot detection}\label{sec:knot-det}
Now, let us discuss how the obtained results affect the issue of distinguishing knots with the help of the colored HOMFLY polynomials. We analyze this question considering specific examples.

\subsection{On $\mathfrak{sl}_2$ kernel}

First, we consider an example of the $\mathfrak{sl}_2$ kernel. This case is more illustrative than the general $\mathfrak{sl}_n$ case, because the kernel first appears at the 4-th order. In addition, in order to detect elements of the $\mathfrak{sl}_2$ kernel, one can use the $\mathfrak{sl}_n$ algebra and the HOMFLY polynomials, which are much better studied and calculated than the Kauffman polynomials needed for the detection of the $\mathfrak{sl}_n$ kernel.

For illustrative purpose, we consider two examples of the $\mathfrak{sl}_2$ kernel in details: order 4 and order 5. The $\mathfrak{sl}_2$ kernel element at the 4-th order appears differently from the $\mathfrak{sl}_n$ case at the 8-th order, which is described in Section \ref{sec:prim-ker-sl}. However, this does not affect the further reasoning and the method of constructing the combination of knots that are detected by this element. Therefore, we analyze this case in as much detail as possible, since the calculations for it are quite simple. At the 5-th level, the $\mathfrak{sl}_2$ kernel element already appears in an absolutely similar way to the $\mathfrak{sl}_n$ kernel element at the 8-th level. We also analyze this case, but omitting some details similar to the ones of the level 4 case.

\subsubsection{Order 4}

An element of the $\mathfrak{sl}_2$ kernel appears at the 4-th order first time. Indeed, $\dim \mathcal{A}_4 = 3$ (see Table \ref{tab:dimensions}) while the number of independent Casimir invariants is 2 (quadratic one and its square). Therefore, there is the 1-dimensional kernel generated by
\be
\label{sl2ker4}
\begin{picture}(150,40)(-20,-20)
\put(0,0){\circle*{40}}
\put(0,0){\color{white}\circle*{37}}
\put(0,0){\circle{20}}
\put(7,7){\line(1,1){7}}
\put(-7,7){\line(-1,1){7}}
\put(7,-7){\line(1,-1){7}}
\put(-7,-7){\line(-1,-1){7}}
\put(7,7){\circle*{2}}
\put(7,-7){\circle*{2}}
\put(-7,7){\circle*{2}}
\put(-7,-7){\circle*{2}}
\put(27,-3){$- \ \, 2$}
\put(70,0){\circle*{40}}
\put(70,0){\color{white}\circle*{37}}
\put(70,0){\line(0,1){20}}
\put(70,0){\line(2,-1){18}}
\put(70,0){\line(-2,-1){18}}
\put(70,0){\circle*{2}}
\put(85,15){\mbox{\fontsize{7.5}{7.5} $2$}}
\put(98,-3){$=\,0$}
\end{picture}
\ee

\

\noindent This relation can be obtained by applying Proposition 6.2 from Vogel's paper \cite{vogel2011algebraic}

\begin{picture}(300,35)(-90,-10)
\put(0,0){$\Phi_{\mathfrak{sl}_2} \Big($}
\put(40,3){\line(-1,1){10}}
\put(40,3){\line(-1,-1){10}}
\put(40,3){\line(1,0){10}}
\put(50,3){\line(1,1){10}}
\put(50,3){\line(1,-1){10}}
\put(65,0){$\Big) = \Phi_{\mathfrak{sl}_2} \ t \Big($}
\put(120,13){\line(1,0){20}}
\put(120,-7){\line(1,0){20}}
\put(145,0){$-$}
\put(160,-7){\line(1,1){20}}
\put(160,13){\line(1,-1){9}}
\put(171,2){\line(1,-1){9}}
\put(185,0){$\Big), \quad \Phi_{\mathfrak{sl}_2}$}
\put(238,3){\circle{25}}
\put(255,0){$= 3$}
\end{picture}

\noindent to the internal graph of the diagram $
\begin{picture}(30,30)(-15,-3)
\put(0,0){\circle*{30}}
\put(0,0){\color{white}\circle*{27}}
\put(0,0){\circle{14}}
\put(5,5){\line(1,1){5}}
\put(-5,5){\line(-1,1){5}}
\put(5,-5){\line(1,-1){5}}
\put(-5,-5){\line(-1,-1){5}}
\put(5,5){\circle*{2}}
\put(5,-5){\circle*{2}}
\put(-5,5){\circle*{2}}
\put(-5,-5){\circle*{2}}
\end{picture}$
and taking into account IHX (see Fig.\,\ref{fig:IHX-relation}) and 1T (Fig.\,\ref{fig:1T}) relations. 

\

Next, we apply STU  relation (Fig.\,\ref{fig:STU-relation}) to closed Jacobi diagrams \eqref{sl2ker4} and convert them into chord diagrams modulo 1T relation:
\be
\label{chrodkersl2}
\begin{picture}(350,40)(-40,-20)
\put(0,0){\circle*{40}}
\put(0,0){\color{white}\circle*{37}}
\put(0,0){\circle{20}}
\put(7,7){\line(1,1){7}}
\put(-7,7){\line(-1,1){7}}
\put(7,-7){\line(1,-1){7}}
\put(-7,-7){\line(-1,-1){7}}
\put(7,7){\circle*{2}}
\put(7,-7){\circle*{2}}
\put(-7,7){\circle*{2}}
\put(-7,-7){\circle*{2}}
\put(27,-3){$- \ \, 2$}
\put(70,0){\circle*{40}}
\put(70,0){\color{white}\circle*{37}}
\put(70,0){\line(0,1){20}}
\put(70,0){\line(2,-1){18}}
\put(70,0){\line(-2,-1){18}}
\put(70,0){\circle*{2}}
\put(85,15){\mbox{\fontsize{7.5}{7.5} $2$}}
\put(98,-3){$=$}
\put(130,0){\circle*{40}}
\put(130,0){\color{white}\circle*{37}}
\put(112,7){\line(1,0){36}}
\put(112,-7){\line(1,0){36}}
\put(123,-18){\line(0,1){36}}
\put(137,-18){\line(0,1){36}}
\put(158,-3){$- \ \ 4$}
\put(200,0){\circle*{40}}
\put(200,0){\color{white}\circle*{37}}
\put(180,0){\line(1,0){40}}
\put(190,-17){\line(0,1){34}}
\put(210,-17){\line(0,1){34}}
\put(200,-20){\line(0,1){40}}
\put(228,-3){$- \ \ 4$}
\put(270,0){\circle*{40}}
\put(270,0){\color{white}\circle*{37}}
\put(250,0){\line(1,0){40}}
\put(270,-20){\line(0,1){40}}
\put(285,15){\mbox{\fontsize{7.5}{7.5} $2$}}

%\put(7,-34){\line(0,1){28}}
%\put(7,-6.5){\circle*{2}}
%\put(18,-4){\mbox{\fontsize{7.5}{7.5} $b$}}
\end{picture}
\ee
Further, we construct singular knots from chord diagrams as described in Section \ref{sec:chord_diag}. The example is given below:
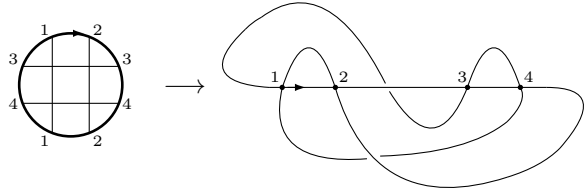
\be
\begin{picture}(50,45)(130,-33)
\put(130,0){\circle*{40}}
\put(130,0){\color{white}\circle*{38}}
\put(112,7){\line(1,0){36}}
\put(112,-7){\line(1,0){36}}
\put(123,-18){\line(0,1){36}}
\put(137,-18){\line(0,1){36}}
\put(128,19.5){\vector(1,0){7}}
\put(165,-3){$\longrightarrow$}
\put(115,19){\mbox{\fontsize{6.5}{6.5} $1$}}
\put(135,19){\mbox{\fontsize{6.5}{6.5} $2$}}
\put(146,8){\mbox{\fontsize{6.5}{6.5} $3$}}
\put(146,-9){\mbox{\fontsize{6.5}{6.5} $4$}}
\put(135,-23){\mbox{\fontsize{6.5}{6.5} $2$}}
\put(115,-23){\mbox{\fontsize{6.5}{6.5} $1$}}
\put(103,-9){\mbox{\fontsize{6.5}{6.5} $4$}}
\put(103,8){\mbox{\fontsize{6.5}{6.5} $3$}}
\end{picture}
\begin{picture}(120,65)(-25,-52)
\put(0,-20){\line(1,0){105}}
\put(5,-20){\circle*{2}}
\put(25,-20){\circle*{2}}
\put(75,-20){\circle*{2}}
\put(95,-20){\circle*{2}}
\qbezier(5,-20)(15,10)(25,-20)
\qbezier(75,-20)(85,10)(95,-20)
\qbezier(0,-20)(-30,-18)(-9,2)
\qbezier(-9,2)(-6,5)(-3,6.5)
\qbezier(-3,6.5)(22,25)(44,-18.5)
\qbezier(45.5,-21.5)(63,-50)(75,-20)
\qbezier(25,-20)(40,-70)(90,-54)
\qbezier(90,-54)(103,-50)(110,-43)
\qbezier(110,-43)(130,-20)(105,-20)
\qbezier(5,-20)(-2,-50)(37,-47)
\qbezier(42,-46)(75,-45)(90,-33)
\qbezier(90,-33)(98,-25)(95,-20)
\put(7,-20){\vector(1,0){7}}
\put(-3,-18){\mbox{\fontsize{6.5}{6.5} $1$}}
\put(23,-18){\mbox{\fontsize{6.5}{6.5} $2$}}
\put(68,-18){\mbox{\fontsize{6.5}{6.5} $3$}}
\put(93,-18){\mbox{\fontsize{6.5}{6.5} $4$}}
\end{picture}
\nonumber
\ee

Note that there are infinitely many singular knots that correspond to each chord diagram. If there are no additional requirements, then we can choose any singular knot that corresponds to the given chord diagram. This is exactly what we did in the above example.

With the help of Vassiliev skein relation (Fig.\,\ref{fig:Vass-skein}) we finally get the following example of linear combination of Vassiliev invariants
\be
\label{vasskersl2}
\mathcal{K}^4_{\mathfrak{sl}_2} := \Big(2 v(4_1) + v(3_1) + v(\bar{3}_1) - 4v(0_1) \Big)  \ - \ 4 \Big( v(5_2) - 3v(3_1) + 3v(0_1) - v(4_1) \Big) \ - \ 4 \Big( v(3_1\texttt{\#} 3_1) - 2v(3_1) + v(0_1) \Big),
\nonumber
\ee
where we use the standard Rolfsen notation for knots, $\bar{3}_1$ is a mirror knot to the trefoil $3_1$, $0_1$ denotes the unknot and $3_1\texttt{\#} 3_1$ is a composite knot of two trefoils. The denotations of knots correspond to pictures from \cite{bar2004knot, knotinfo}. %Pictures of corresponding knots, which we imply here to distinguish them, can be found in \cite{bar2004knot, knotinfo}.

If we consider the Vassiliev invariant of order 2 ($v \in \mathcal{V}_2$) or 3 ($v \in \mathcal{V}_3$), then this combination \textit{vanishes} $\mathcal{K}^4_{\mathfrak{sl}_2} = 0$, which follows from the definition of Vassiliev invariants \ref{def:Vass_inv_def}. 

If we consider a Vassiliev invariant of order 4 ($v \in \mathcal{V}_4$), which is obtained from the (colored) Jones polynomials (that correspond to $\mathfrak{sl}_2$ weight system), then this combination \textit{also vanishes}. It follows from two facts. First, a value of $v \in \mathcal{V}_n$ on a knot $\mathcal{K}$ with $n$ singular points depends only on the chord diagram of $\mathcal{K}$. Second, this singular knot corresponds to the combination of chord diagrams \eqref{chrodkersl2} that belongs to the kernel of $\mathfrak{sl}_2$ weight system. %$\mathcal{K}^4_{\mathfrak{sl}_2}$ corresponds to combination of chord diagrams \eqref{chrodkersl2} belonging to the kernel of $\mathfrak{sl}_2$-weight system. 
For another validation, let us calculate the perturbative expansion of the HOMFLY polynomials (in fundamental representation) for these knots:
\be
\label{jonesexample}
H^{\mathcal{K}^4_{\mathfrak{sl}_2}} = -48(n-1)(n-2) (5n+18) (n+3) \hbar^4 \ -  \ 384n(n-1)(n+3)(11n^2+40n-2)\hbar^5 \ + \ O(\hbar^6)\,.
\ee
Recall that the HOMFLY polynomial corresponds to $\mathfrak{sl}_n$ weight system (see Section~\ref{sec:CS-KI}). This formula illustrates that the linear combination $\mathcal{K}^4_{\mathfrak{sl}_2}$ cannot be detected by Vassiliev invariants of order $\leq 4$ coming from $\mathfrak{sl}_2$ quantum knot invariants (the Jones polynomials), while Vassiliev invariants of order $> 4$ do detect this combination already for $\mathfrak{sl}_2$ weight system (note that the 5-th order in \eqref{jonesexample} is not zero for $n=2$). 

Thus, the existence of a kernel of a weight system does not itself guarantee that one is unable to distinguish some knots or their combinations via higher order Vassiliev invariants or entire quantum knot invariants. The fact that there are knots that are not distinguishable by $\mathfrak{sl}_2$ invariants (e.g. mutant knots) requires additional considerations.

\subsubsection{Order 5}

At the 5-th order, the element of the $\mathfrak{sl}_2$ weight system kernel can be constructed in the same way as for the kernel of $\mathfrak{sl}_n$, what was described in Section \ref{method}. Indeed, let us consider the diagram $D_5 := \frac14 \hat{\omega}\bbl$. This diagram is a primitive one of the 5-th order. According to Proposition 6.2 from Vogel's paper \cite{vogel2011algebraic} $\omega=0$ for the algebra $\mathfrak{sl}_2$, therefore, $D_5$ vanishes for $\mathfrak{sl}_2$. To draw this diagram, we recall (Remark \ref{remark32}) that $\hat{\omega} = \frac{8}{3}\hat{t}^3 - \frac{2}{3}\hat{x}_3$, so
\begin{equation}
	\begin{tikzpicture}[scale=0.7]
		\coordinate (a) at (0, 0);
		\coordinate (b) at (3.5, 0);
		\node at ($0.5*(a)+0.5*(b)$) {\Large $-\frac13$};
		\node at ($(a)+(-2.5,0)$) {\Large $D_5 = \frac16$};
		
		\draw[thick] (a) circle(1);
		\bubble{(a)}{($(a)+(120:1)$)};
		\bubble{(a)}{($(a)+(240:1)$)};
		\bubble{(a)}{($(a)+(0:1)$)};
		
		\draw[thick] (b) circle(1);
		\xVogel{($(b)+(120:1)$)}{($(b)+(240:1)$)}{($(b)+(1,0)$)}{3};
	\end{tikzpicture}
\end{equation}

\

 Next, we construct a singular knot in the same way as described in the previous subsection. For brevity, we omit details here and write down the answer, i.e. a specific example of a linear combination of Vassiliev invariants constructed according to the $D_5$ diagram:
\be
\mathcal{K}^5_{\mathfrak{sl}_2} := 11v(0_1) + 12v(3_1\texttt{\#}\bar{3}_1) + 5v(4_1\texttt{\#}\bar{3}_1) - 8v(3_1) - 4v(4_1) + 5v(5_1) - 10v(5_2) - 5v(6_1) + 5v(6_2) + 7v(7_6) -  \nonumber \\
-20v(\bar{3}_1) +v(\bar{5}_2) + v(\bar{6}_1) + 7v(\bar{6}_3) - 7v(\bar{8}_{20})\,. \hspace{3mm} 
\ee
To find the values of Vassiliev invariants let us calculate the HOMFLY polynomials taken in the fundamental representation for these knots:
\be
H^{\mathcal{K}^5_{\mathfrak{sl}_2}} = -n(n-1)(n-2)(n+2)(n+1)\hbar^5 - \frac{1}{2}n^2(n-1)(n+1)(13n^2-4)\hbar^6 + O(\hbar^8)\,.
\ee
We see that for $n=2$ the 5-th order vanishes while the higher orders do not.

\subsection{On $\mathfrak{sl}_n$ kernel}

All the above reasoning is also completely true for the $\mathfrak{sl}_n$ kernel. Combinations of knots that can be constructed using the kernel element $\widehat{(\omega P_{sl})}\bbl$ would not differ by Vassiliev invariants of orders $\leq 8$ coming from the HOMFLY polynomial. However, for orders $\geq 9$ this would no longer be the case.

How to pass from some linear combination of knots to just a difference between two knots $K - K'$? We do not know an algorithmic answer to this question. However, from a chord diagram of the $n$-th order, one can construct infinitely many singular knots with $n$ double points. Further, these singular knots can be combined with each other (added and subtracted with arbitrary coefficients). It may well happen that as a result of these manipulations, we get a linear combination of knots of the form $K - K'$ for some knots $K$ and $K'$. For small orders and knots crossing numbers, this can be done manually.

To summarize, we can say that both the colored HOMFLY polynomials and Vassiliev invariants might simultaneously be complete knot invariants. And they might be expressed through each other in a very complex way. At least today, this does not contradict anything including the result of P. Vogel.

\section*{Acknowledgements}

We are grateful for enlightening discussions to A. Morozov and N. Tselousov. 

This work was funded by the RSF grant No.24-12-00178.

\printbibliography
 
\newpage
\appendixpage

\begin{appendices}
	\section{Vogel algebra: advanced features} \label{appendix}

	In this appendix, more sophisticated attributes of Vogel algebra are given. They all were in Vogel's paper \cite{vogel2011algebraic}, but many of them were not written explicitly. Among those are the \textit{marked diagrams} that stand for some linear combinations of regular trivalent diagrams modulo the AS and IHX relations. They contain a 4-valent vertex marked by some non-negative integer number $p \in \mathbb{N}_0$. Marked diagrams are defined as follows:
	
	\begin{figure}[h]
		\centering
		\begin{tikzpicture}[scale=0.5]
			\coordinate (a) at (0,0); 
			\coordinate (b) at (3.5, 0);
			\coordinate (c) at (7, 0);
			\coordinate (d) at (10.5, 0);
			\coordinate (e) at (14, 0);
			
			\node at ($0.5*(a)+0.5*(b)$) {$=2\hat{t}^2$};
			\node at ($0.5*(b)+0.5*(c)$) {$+2\hat{t}$};
			\node at ($0.5*(c)+0.5*(d)$) {$-\frac43$};
			\node at ($0.5*(d)+0.5*(e)$) {$-\frac23$};
			
			\draw (a)+ (-1, 1)-- +(1, 1);
			\draw (a) +(-1, -1)-- +(1, -1);
			\draw (a) +(0, 1)-- +(0, -1);
			\draw[fill] (a) circle (0.1);
			\node[scale=0.7] at ($(a)+(0.3, 0)$) {0};
			
			\draw (b)+ (-1, 1)-- +(1, 1);
			\draw (b) +(-1, -1)-- +(1, -1);
			\draw (b) +(0, 1)-- +(0, -1);
			
			\draw (c)+ (-1, 1)-- +(1, 1);
			\draw (c) +(-1, -1)-- +(1, -1);
			\draw (c) +(0.3, 1)-- +(0.3, -1);
			\draw (c) +(-0.3, 1)-- +(-0.3, -1);

			\draw (d) +(-1, 1) -- + (1, 1);
			\draw (d) +(-1, -1) -- + (1, -1);
			\draw (d) +(0.5, -1) -- + (0.5, 1);
			\draw (d) +(-0.5, -1) -- + (-0.5, 1);
			\draw (d) +(0, -1) -- + (0, 1);
			
			\draw (e)+ (-1, 1)-- +(1, 1);
			\draw (e) +(-1, -1)-- +(1, -1);
			\draw (e) +(0.3, 1)-- +(0.3, -1);
			\draw (e) +(-0.3, 1)-- +(-0.3, -1);
			\draw (e) +(-0.3, 0) -- +(0.3, 0);
		\end{tikzpicture} \\
		\vspace{0.4cm}
		\begin{tikzpicture}[scale=0.5]
			\coordinate (a) at (-14.5, 0);
			\coordinate (b) at (-8.5, 0);
			\coordinate (c) at (-2, 0);
			\coordinate (d) at (2, 0);
			\coordinate (e) at (6, 0);
			\coordinate (f) at (10, 0);
			\coordinate (g) at (14, 0);
			\node at ($(a) + (-2, 0)$) {$\varepsilon_0 =$};
			\node at ($0.5*(a) +0.5*(b)$) {$= \left( \hat{\omega} -\frac{8 \hat{t}^3}{3} \right) $};
			\node at ($0.5*(b) +0.5*(c)$) {$- \left( \hat{\omega} -\frac{4 \hat{t}^3}{3} \right) $};
			\node at ($0.5*(c) +0.5*(d)$) {$+ \frac{2 \hat{t}^2}{3} $};
			\node at ($0.5*(d) +0.5*(e)$) {$+ \frac{10 \hat{t}}{3} $};
			\node at ($0.5*(e) +0.5*(f)$) {$- \frac{4}{3} $};
			\node at ($0.5*(f) +0.5*(g)$) {$- \frac{2}{3}$};
			%	 epsilon0
			\draw[fill] (a) circle (0.1);
			\draw (a) -- + (1,-1);
			\draw (a) -- + (1,1);
			\draw (a) -- + (-1,1);
			\draw (a) -- + (-1,-1);
			\node[scale=0.7] at ($(a)+(0, -0.5)$) {$0$};
			
			% PSI
			\draw (b) +(-1, 1) -- + (1, 1);
			\draw (b) +(-1, -1) -- + (1, -1);
			\draw (b) +(0, -1) -- + (0, 1);
			
			% H
			\draw (c) +(-1, -1) -- + (-1, 1);
			\draw (c) +(1, -1) -- + (1, 1);
			\draw (c) +(-1, 0) -- + (1, 0);
			
			% PSI2
			\draw (d) +(-1, 1) -- + (1, 1);
			\draw (d) +(-1, -1) -- + (1, -1);
			\draw (d) +(0.3, -1) -- + (0.3, 1);
			\draw (d) +(-0.3, -1) -- + (-0.3, 1);
			
			% PSI3
			\draw (e) +(-1, 1) -- + (1, 1);
			\draw (e) +(-1, -1) -- + (1, -1);
			\draw (e) +(0.5, -1) -- + (0.5, 1);
			\draw (e) +(-0.5, -1) -- + (-0.5, 1);
			\draw (e) +(0, -1) -- + (0, 1);
			
			% PSI4
			\draw (f) +(-1, 1) -- + (1, 1);
			\draw (f) +(-1, -1) -- + (1, -1);
			\draw (f) +(0.6, -1) -- + (0.6, 1);
			\draw (f) +(-0.6, -1) -- + (-0.6, 1);
			\draw (f) +(0.2, -1) -- + (0.2, 1);
			\draw (f) +(-0.2, -1) -- + (-0.2, 1);
			
			% wierd	
			
			\draw (g) +(-1, 1) -- + (1, 1);
			\draw (g) +(-1, -1) -- + (1, -1);
			\draw (g) +(0.5, -1) -- + (0.5, 1);
			\draw (g) +(-0.5, -1) -- + (-0.5, 1);
			\draw (g) +(0, -1) -- + (0, 1);
			\draw (g) +(0, 0) -- + (0.5, 0);
			
		\end{tikzpicture}
		\caption{Definition of diagrams marked by zero.}
		\label{fig:marked0}
	\end{figure}
	
		\begin{figure}[h!]
		\centering
		\begin{tikzpicture}[scale=0.5]
			\coordinate (a) at (-0.5,0); 
			\coordinate (b) at (3.5, 0);
			\coordinate (c) at (7, 0);
			\coordinate (d) at (10.5, 0);
			\coordinate (e) at (14, 0);
			
			\node at ($0.5*(a)+0.5*(b)$) {$=\frac23$};
			\node at ($0.5*(b)+0.5*(c)$) {$-\frac23$};
			\node at ($0.5*(c)+0.5*(d)$) {$+\frac{8 \hat t^2}{9}$};
			\node at ($0.5*(d)+0.5*(e)$) {$+\frac{4\hat t \hat\omega_p}{3}$};
			
			\draw (a)+ (-1, 1)-- +(1, 1);
			\draw (a) +(-1, -1)-- +(1, -1);
			\draw (a) +(0, 1)-- +(0, -1);
			\draw[fill] (a) circle (0.1);
			\node[scale=0.7] at ($(a)+(0.7, 0)$) {$p{+}1$};
			
			\draw (b)+ (-1, -1)-- +(1, 1);
			\draw (b) +(1, -1)-- +(-1, 1);
			\draw[fill] (b) circle(0.1);
			\node[scale=0.7] at ($(b)+(0, -0.4)$) {$p$};
			\draw (b)+(45:0.5) arc(45:130:0.5);
			\draw (b)+(140:0.5) arc(140:225:0.5);

			\draw (c)+ (-1, -1)-- +(1, 1);
			\draw (c) +(1, -1)-- +(-1, 1);
			\draw[fill] (c) circle(0.1);
			\node[scale=0.7] at ($(c)+(0, -0.4)$) {$p$};
			\draw (c)+(145:0.5) arc(145:225:0.5);

			\draw (d)+ (-1, 1)-- +(1, 1);
			\draw (d) +(-1, -1)-- +(1, -1);
			\draw (d) +(0, 1)-- +(0, -1);
			\draw[fill] (d) circle (0.1);
			\node[scale=0.7] at ($(d)+(0.4, 0)$) {$p$};
			
			\draw (e)+ (-1, 1)-- +(1, 1);
			\draw (e) +(-1, -1)-- +(1, -1);
			\draw (e) +(0, 1)-- +(0, -1);
		\end{tikzpicture}
		%	\vspace{1cm}
		\begin{picture}(750,30)(-80, -20)
			\put(-28,-17){\mbox{$\varepsilon_p \ =$}}
			
			\put(0,-30){\line(1,1){30}}
			\put(0,0){\line(1,-1){30}}
			\put(15,-15){\circle*{4}}
			\put(13,-25){$p$}
			\put(43,-17){=}
			\put(73,-30){\line(0,1){30}}
			\put(83,-30){\line(0,1){30}}
			\put(63,0){\line(1,0){30}}
			\put(63,-30){\line(1,0){30}}
			\put(83,-15){\circle*{4}}
			\put(85,-25){$p$}
			\put(100,-17){$- \ \hat\omega_p$}
			\put(128,-30){\line(0,1){30}}
			\put(158,-30){\line(0,1){30}}
			\put(128,-15){\line(1,0){30}}
			\put(165,-17){$+ \ \hat\omega_p$}
			\put(208,-30){\line(0,1){30}}
			\put(193,0){\line(1,0){30}}
			\put(193,-30){\line(1,0){30}}
			\put(228,-17){$- \ \frac{4}{3}\hat t$}
			\put(271,-30){\line(0,1){30}}
			\put(256,0){\line(1,0){30}}
			\put(256,-30){\line(1,0){30}}
			\put(271,-15){\circle*{4}}
			\put(273,-25){$p$}
			\put(291,-17){$+ \ \frac{2}{3}\hat t$}
			\put(319,-30){\line(0,1){30}}
			\put(349,-30){\line(0,1){30}}
			\put(319,-15){\line(1,0){30}}
			\put(334,-15){\circle*{4}}
			\put(336,-25){$p$}
		\end{picture}
		\vspace{0.15cm}
		\caption{Definition of diagrams marked by $p>0$.}
		\label{fig:marked_recursion}
	\end{figure}

	\noindent The diagram $\hat\omega_p$ in terms of these marked diagrams can be displayed as follows:

	\begin{figure}[h!]
		\centering
		\begin{tikzpicture}[scale=0.5]
			\draw (90:1)--(90:2);
			\draw (210:1)--(210:2);
			\draw (330:1)--(330:2);
			\draw (90:1) -- (210:1) -- (330:1) -- (90:1);
			\draw [fill] (0, -0.5) circle(0.1);
			\node [scale=0.7] at (0, -0.9) {$p$};
			
		\end{tikzpicture}
		\caption{$\hat{\omega}_p$ as a marked diagram.}
		\label{fig:omega_p}
	\end{figure}

	These marked diagrams have an important property. Namely, a connected diagram having more than one marked point depends only on its total marking. This is why the relation $\hat{ \omega}_p \hat{ \omega}_q$ = $\hat{ \omega}_0 \hat{ \omega}_{p+q}$ holds. And this is also a reason for  the multiplication by $\sigma$ to make sense as an operation that raises the total marking by $1$.
	
	The diagram $\varepsilon_p$ defined above is symmetric with respect to the permutations of legs. Hence, it corresponds to itself as an open Jacobi diagram.

	For a Lie algebra weight system there are two additional relations that involve marked diagrams. They provide a way to multiply diagrams by $\sigma$ but only under the action of a Lie algebra weight system, see Fig.\,\ref{fig:two_additional_relations}.

	\begin{figure}[h!]
		\centering
		\begin{tikzpicture}[scale=0.56]
			\coordinate (a) at (-9, 0);
			\coordinate (b) at ($(a)+(3.25, 0)$);
			\coordinate (c) at ($(b)+(2, 0)$);
			\coordinate (d) at ($(c)+(4.5, 0)$);
			\coordinate (e) at ($(d)+(3.25, 0)$);
			\coordinate (f) at ($(e)+(2, 0)$);

			\node at ($(a)+(-2.25, 0)$) { $\sigma_L \Phi_L$};
			\node at ($(a)+(-1.1, 0)$) {\Huge $($};
			\node at ($0.55*(a)+0.45*(b)$) {$-t_L$};
			\node at ($(b)+(-0.95, 0)$) {\huge $(\,$};
			\node at ($0.5*(b)+0.5*(c)$) {$-$};
			\node at ($(c)+(0.95, 0)$) {\huge $\,)$};
			\node at ($(c)+(1.1, 0)$) {\Huge $ \, \ ) \ $};
			\node at ($0.5*(c)+0.5*(d)$) {$=\,\Phi_L$};
			\node at ($(d)+(-1.1, 0)$) {\Huge $($};
			\node at ($0.55*(d)+0.45*(e)$) {$-\omega_L$};
			\node at ($(e)+(-0.95, 0)$) {\huge $(\,$};
			\node at ($0.5*(e)+0.5*(f)$) {$-$};
			\node at ($(f)+(0.95, 0)$) {\huge $\,)$};
			\node at ($(f)+(1.1, 0)$) {\Huge $\,)$};

			\draw (a)+(45:-1)--+(-0.4, 0)--+(135:1);
			\draw (a)+(-0.4, 0) -- +(0.4, 0);
			\draw (a)+(45:1) -- +(0.4, 0) -- +(135:-1);
			
			\draw (b)+ (135:1)-- +(45:1);
			\draw (b) +(135:-1)-- +(45:-1);
			
			\draw (c)+(45:1) -- +(45:-1);
			\draw (c)+(135:-1)-- +(135:-0.1);
			\draw (c) +(135:0.1) -- +(135:1);
			
			\draw (d)+(45:-1)--+(-0.4, 0)--+(135:1);
			\draw (d)+(-0.4, 0) -- +(0.4, 0);
			\draw (d)+(45:1) -- +(0.4, 0) -- +(135:-1);
			\draw [fill] (d) circle(0.1);
			\node[scale=0.7] at ($(d)+(0, -0.3)$) {$0$};
			
			\draw (e)+ (135:1)-- +(45:1);
			\draw (e) +(135:-1)-- +(45:-1);
			
			\draw (f)+(45:1) -- +(45:-1);
			\draw (f)+(135:-1)-- +(135:-0.1);
			\draw (f) +(135:0.1) -- +(135:1);

		\end{tikzpicture}\\
		\vspace{0.4cm}
		\begin{tikzpicture}[scale=0.56]
			\coordinate (a) at (-9.5, 0);
			\coordinate (b) at ($(a)+(3,0)$);
			\coordinate (c) at ($(b)+(2,0)$);
			\coordinate (d) at ($(c)+(3.75,0)$);
			\coordinate (e) at ($(d)+(2.,0)$);
			\coordinate (f) at ($(e)+(2.,0)$);
			\coordinate (k) at ($(f)+(4.25,0)$);
			\coordinate (l) at ($(k)+(4,0)$);
			\coordinate (m) at ($(l)+(2.,0)$);
			\coordinate (n) at ($(m)+(2.,0)$);

			\node at ($(a)+(-2.25, 0)$) {$\sigma_L \Phi_L$}; 
			\node at ($(a)+(-1.10,0)$) {\Huge $($};
			\node at ($0.55*(a)+0.45*(b)$) {$-\frac{t_L}{3}\, $};
			\node at ($(b)+(-0.95, 0)$) {\huge $(\,$};
			\node at ($0.5*(c)+0.5*(b)$) {$+$};
			
			\node at ($(c)+(0.95,0) $) {\huge $\,)\,$};
			\node at ($0.5*(c)+0.5*(d)$) {$-\frac{2t_L^2}{3}\, $};
			\node at ($(d)+(-0.95, 0)$) {\huge $(\,$};
			
			\node at ($0.5*(e)+0.5*(d)$) {$+$};
			\node at ($0.5*(e)+0.5*(f)$) {$+$};
			
			\node  at ($(f)+(0.95, 0)$) {\huge $\,)$};
			\node  at ($(f)+(1.1, 0)$) {\Huge $\, \ ) \ $};
			\node at ($0.5*(k)+0.5*(f)$) {$ \, = \Phi_L$};
			\node at ($(k)+(-1.1, 0)$) {\Huge $($};
			
			\node at ($0.55*(k)+0.45*(l)$) {$-\frac{2t_L \omega_L}{3}\, $};
			\node at ($(l)+(-0.95,0)$) {\huge $(\,$};
			
			\node at ($0.5*(m)+0.5*(l)$) {$+$};
			\node at ($0.5*(m)+0.5*(n)$) {$+$};
			\node at ($(n)+(0.95,0)$) {\huge $\,)$};
			\node at ($(n)+(1.1,0)$) {\Huge $\,)$};

			\draw (a) circle (0.4);
			\draw (a) +(45:0.4) -- +(45:1);
			\draw (a) +(135:0.4) -- +(135:1);
			\draw (a) +(225:0.4) -- +(225:1);
			\draw (a) +(315:0.4) -- +(315:1);
			
			\draw (b) +(45:1)-- +(0, 0.4);
			\draw (b) +(135:1)-- +(0, 0.4);
			\draw (b) +(45:-1)-- +(0, -0.4);
			\draw (b) +(135:-1)-- +(0, -0.4);
			\draw (b) +(0, -0.4)-- +(0, 0.4);
			
			\draw (c) +(45:1)-- +( 0.4,0);
			\draw (c) +(135:1)-- +(-0.4, 0);
			\draw (c) +(45:-1)-- +(-0.4,0);
			\draw (c) +(135:-1)-- +(0.4,0);
			\draw (c) +(-0.4,0)-- +(0.4, 0);
			
			\draw (d) +(225:1)-- +(135:1);
			\draw (d) +(315:1) -- +(45:1);
			
			\draw (e) +(135:1) -- +(45:1);
			\draw (e) +(225:1) -- +(315:1);
			
			\draw (f) +(45:1)-- + (45:-1);
			\draw (f) +(135:0.1) -- +(135:1);
			\draw (f) +(135:-0.1) -- +(135:-1);

			\draw (l) +(225:1)-- +(135:1);
			\draw (l) +(315:1) -- +(45:1);
			
			\draw (m) +(135:1) -- +(45:1);
			\draw (m) +(225:1) -- +(315:1);
			
			\draw (n) +(45:1)-- + (45:-1);
			\draw (n) +(135:0.1) -- +(135:1);
			\draw (n) +(135:-0.1) -- +(135:-1);
			
			\draw (k) +(45:1)-- + (45:-1);
			\draw (k) +(135:-1) -- +(135:1);
			\draw [fill] (k) circle (0.07);
			\node[scale =0.6] at ($(k)+(0,-0.3)$) {$0$};

		\end{tikzpicture}
		\caption{Additional relations that work for a Lie algebra weight systems.}
		\label{fig:two_additional_relations}
	\end{figure}
	
\end{appendices}

\end{document}